\documentclass{article}
\usepackage{amsthm,setspace}
\usepackage{epsfig,amsmath,amssymb}
\usepackage[all]{xy}
\usepackage{lscape}
\usepackage{graphicx}

\usepackage [UglyObsolete,tight,heads=LaTeX] {diagrams}
\newcommand{\dd}{\dLine}
\newcommand{\rd}{\rdLine}
\newcommand{\ld}{\ldLine}

\newtheorem{teo}{Theorem}[section]

\newtheorem{cor}[teo]{Corollary}
\newtheorem{pro}[teo]{Proposition}
\newtheorem{lem}[teo]{Lemma}

\theoremstyle{definition}
\newtheorem{defi}[teo]{Definition}

\title{Right coideal subalgebras of quantized universal enveloping algebras of type $G_2$\thanks{This paper was written based on a research period at the UNAM FES-C, Mexico, with the support of CNPq-Brazil. It is part of the author's PhD thesis.}}

\author{B\'arbara Pogorelsky\\\footnotesize{barbara.pogorelsky@ufrgs.br}\\\footnotesize{Instituto de Matem\'atica, Universidade Federal do Rio Grande do Sul}\\\footnotesize{Av. Bento Gon\c calves 9500, Porto Alegre, RS, 91509-900, Brazil}}

\begin{document}

\maketitle

\begin{abstract}
\noindent In this paper we describe the right coideal subalgebras
containing all group-like elements of the two-parameter quantum
group $U_q(\mathfrak{g})$, where $\mathfrak{g}$ is a simple Lie
algebra of type $G_2$, while the main parameter of quantization
$q$ is not a root of 1. As a consequence, we determine that there are precisely 60 different right coideal
subalgebras containing all group-like elements. If the multiplicative order $t$ of $q$ is
finite, $t>4$, $t\neq 6$, then the same classification remains
valid for homogeneous right coideal subalgebras of the two-parameter version of the
small Lusztig quantum group $u_q(\mathfrak{g})$.
\end{abstract}



\section{Introduction}

\quad In this paper we continue the classification of the right coideal subalgebras of quantized enveloping algebras started in \cite{Kh6}, \cite{Kh8} and \cite{BSP} by V. K. Kharchenko, A. V. Lara Sagah\'on and the author. In the second section, following \cite{Kh6} and \cite{BSP}, we introduce the main concepts and general results that are necessary.

In the third section we present the algebra $U_q(\mathfrak{g})$, where $\mathfrak{g}$ is a simple Lie
algebra of type $G_2$, which has a triangular decomposition $$U_q(\mathfrak{g})=U_q^-(\mathfrak{g})\otimes_{\textbf{k}[F]}\textbf{k}[H]\otimes_{\textbf{k}[G]}U_q^+(\mathfrak{g}).$$ Here we also give some previous important results obtained in \cite{BSP} for the algebra $U_q^+(\mathfrak{g})$. In particular, Theorem \ref{pbw1} calculates a PBW-basis for this algebra, and Theorem \ref{lattice} completely describes the right coideal subalgebras containing all group-like elements, proving a conjecture of V. K. Kharchenko (see \cite{Kh8}) for the specific case where $\mathfrak{g}$ is a simple Lie algebra of type $G_2$. Recently, this conjecture was proved in general by I. Heckenberger and H.-J. Schneider in \cite{HS}. In this section we also prove in Lemma \ref{lema1} that the right coideal subalgebras that contain all group-like elements have a triangular decomposition $$\textbf{U}=\textbf{U}^-\otimes_{\textbf{k}[F]}\textbf{k}[H]\otimes_{\textbf{k}[G]}\textbf{U}^+$$ and in Lemma \ref{lema2} we present necessary and sufficient conditions to determine when the above tensor product of right coideal subalgebras $\textbf{U}^-$ and
$\textbf{U}^+$ of, respectively, $U_q^-(\mathfrak{g})$ and $U_q^+(\mathfrak{g})$ is a right coideal subalgebra of $U_q(\mathfrak{g})$.

In the fourth section we have Theorem \ref{teorema}, which is the main result of the paper. Here we prove, using Lemmas \ref{lema1} and \ref{lema2} that we have the following right coideal subalgebras for $U_q(\mathfrak{g})$

\vspace{0.5cm}
\begin{tabular}{|c|c|c|c|c|c|c|c|c|c|c|c|c|}
  \hline
   & $U_{11}$ & $U_{10}$ & $U_9$ & $U_8$ & $U_7$ & $U_1$ & $U_2$ & $U_3$ & $U_4$ & $U_5$ & $U_6$ & $U_{12}$ \\
   \hline
  $U^-_{11}$ &   &   &   &   &   & $\checkmark$ &   &   &  &   &   &   \\
  \hline
  $U^-_{10}$ &   &   &   &   &   & $\checkmark$ & $\checkmark$  & $\checkmark$  &   &   &   &   \\
  \hline
  $U^-_9$ &   &   &   &   &   & $\checkmark$ & $\checkmark$  & $\checkmark$  & $\checkmark$  &   &   &   \\
  \hline
  $U^-_8$ &   &   &   &   &   & $\checkmark$ & $\checkmark$  & $\checkmark$  & $\checkmark$  & $\checkmark$  & $\checkmark$  &   \\
  \hline
  $U^-_7$ &   &   &   &   & $\checkmark$  & $\checkmark$ & $\checkmark$ & $\checkmark$  & $\checkmark$ & $\checkmark$ & $\checkmark$ & $\checkmark$\\
  \hline
  $U^-_1$ & $\checkmark$ & $\checkmark$ & $\checkmark$ & $\checkmark$ & $\checkmark$ & $\checkmark$ & $\checkmark$ & $\checkmark$ & $\checkmark$ & $\checkmark$ & $\checkmark$ & $\checkmark$ \\
  \hline
  $U^-_2$ &   & $\checkmark$ & $\checkmark$ & $\checkmark$ & $\checkmark$ & $\checkmark$ & $\checkmark$ &   &   &   &   & $\checkmark$ \\
  \hline
  $U^-_3$ &   & $\checkmark$ & $\checkmark$ & $\checkmark$ & $\checkmark$ & $\checkmark$ &   &   &   &   &   &   \\
  \hline
  $U^-_4$ &   &   & $\checkmark$ & $\checkmark$ & $\checkmark$ & $\checkmark$ &   &   &   &   &   &   \\
  \hline
  $U^-_5$ &   &   &   & $\checkmark$ & $\checkmark$ & $\checkmark$ &   &   &   &   &   &   \\
  \hline
  $U^-_6$ &   &   &   & $\checkmark$ & $\checkmark$ & $\checkmark$ &   &   &   &   &   &   \\
  \hline
  $U^-_{12}$ &   &   &   &   & $\checkmark$ & $\checkmark$ & $\checkmark$ &   &   &   &   & $\checkmark$ \\
  \hline
\end{tabular}

\vspace{0.5cm}
\noindent where $\textbf{U}^+=U_i$ and $\textbf{U}^-=U_j^-,$ $1\leq i,j\leq 12$,
are the right coideal subalgebras of $U_q^+(\mathfrak{g})$ and $U_q^-(\mathfrak{g})$ described in \cite{BSP} (see Figure 1 and notations on pages 12-13). In particular, we see that $U_q
({\mathfrak g})$ has precisely 60 different right coideal subalgebras that
contain all group-like elements. In the case where $q$ is a root of $1$ of order $t,$ $t>4,$ $t\neq 6$, the same results are valid for the homogeneous right coideal subalgebras of $u_q(\mathfrak{g})$.

\section{Preliminaries}

\quad In this section, following \cite{Kh6} and \cite{BSP}, we recall the basic concepts.

\begin{defi}
Let $S$ be an algebra over a field $\textbf{k}$ and $A$ its
subalgebra with a fixed basis $\{a_j|j\in J\}$. A linearly ordered
subset $W \subseteq S$ is said to be a set of
\textit{PBW-generators of $S$ over $A$} if there exists a function
$h:W\rightarrow\mathbb{Z}^+\cup\infty$, called the height
function, such that the set of all products
\begin{equation}\label{pbw}
a_jw_1^{n_1}w_2^{n_2}\ldots w_k^{n_k},
\end{equation}
where $j \in J$, $w_1<w_2<\ldots <w_k \in W$, $n_i<h(w_i)$, $1\leq
i\leq k$ is a basis of $S$. The value $h(w)$ is referred to as the
\textit{height} of $w$ in $W$. If $A=\textbf{k}$ is the ground
field, then we shall call $W$ simply as a set of PBW-generators of
$S$.
\end{defi}

\begin{defi}
Let $W$ be a set of PBW-generators of $S$ over a subalgebra $A$.
Suppose that the set of all words in $W$ as a free monoid has its
own order $\prec$ (that is, $a \prec b$ implies $cad \prec cbd$
for all words $a,b,c,d \in W$). A \textit{leading word} of $s \in
S$ is the maximal word $m=w_1^{n_1}w_2^{n_2}\ldots w_k^{n_k}$ that
appears in the decomposition of $s$ in the basis (\ref{pbw}). A
\textit{leading term} of $s$ is the sum $am$ of all terms
$\alpha_ia_im$ that appear in the decomposition of $s$ in the
basis (\ref{pbw}), where $m$ is the leading word of $s$.
\end{defi}

\begin{defi}
A Hopf algebra $H$ is said a \textit{character Hopf algebra} if
the group $G$ of all group-like elements is commutative and $H$ is
generated over $\textbf{k}[G]$ by skew primitive semi-invariants
$a_i, i \in I$: $$\Delta(a_i)=a_i\otimes 1+g_i\otimes a_i,\quad
g^{-1}a_ig=\chi^i(g)a_i,\quad g, g_i \in G,$$ where $\chi^i$, $i
\in I$, are characters of the group $G$.
\end{defi}

Let us associate a quantum variable $x_i$ to $a_i$. For each word
$u$ in $X=\{x_i|i\in I\}$ we denote by $g_u$ an element of $G$
that appears from $u$ by replacing each $x_i$ with $g_i$. In the
same way we denote by $\chi^u$ a character that appears from $u$
by replacing each $x_i$ with $\chi^i$. We define a bilinear skew
commutator on homogeneous linear combinations of words by the
formula
\begin{equation}\label{brac}
[u,v]=uv-\chi^u(g_v)vu,
\end{equation}
where we use the notation $\chi^u(g_v)=p_{uv}=p(u,v)$.

The group $G$ acts on the free algebra $\textbf{k}\langle X\rangle$ by
$g^{-1}ug=\chi(g)u$, where $u$ is an arbitrary monomial in $X$.
The skew group algebra $G\langle X\rangle$ has the natural Hopf
algebra structure
$$\Delta(x_i)=x_i\otimes 1+g_i\otimes x_i,\quad i \in I,\quad
\Delta(g)=g\otimes g.$$ We fix a Hopf algebra homomorphism
$$\xi :G\langle X\rangle\rightarrow H,\quad \xi(x_i)=a_i,\quad
\xi(g)=g,\quad i\in I,\quad g\in G.$$

\begin{defi}
A \textit{constitution} of a word $u$ in $G\cup X$ is a family of
non-negative integers $\{m_x,x\in X\}$ such that $u$ has $m_x$
occurrences of $x$. Certainly almost all $m_x$ in the constitution
are zero.
\end{defi}

Let us fix an arbitrary complete order $<$ on the set $X$, and let
$\Gamma^+$ be the free additive (commutative) monoid generated by
$X$. The monoid $\Gamma^+$ is a completely ordered monoid with
respect to the following order:
\begin{equation}\label{ordemgrau}
m_1x_{i_1}+m_2x_{i_2}+\ldots
+m_kx_{i_k}>m_1'x_{i_1}+m_2'x_{i_2}+\ldots +m_k'x_{i_k}
\end{equation}
if the first from the left nonzero number in
$(m_1-m_1',m_2-m_2',\ldots,m_k-m_k')$ is positive, where
$x_{i_1}>x_{i_2}>\ldots>x_{i_k}$ in $X$. We associate a formal
degree $D(u)=\sum_{x \in X}m_xx \in \Gamma^+$ to a word $u$ in
$G\cup X$, where $\{m_x|x\in X\}$ is the constitution of $u$.
Respectively, if $f=\sum \alpha_iu_i \in G\langle X\rangle$,
$0\neq \alpha_i \in \textbf{k}$ then
\begin{equation}\label{grau}
D(f)=max_i\{D(u_i)\}. \end{equation}
On the set of all words in $X$
we fix the lexicographical order with the priority from the left
to the right, where a proper beginning of a word is considered to
be greater than the word itself.

According to \cite[Theorem 1.1]{Kh5}, every right coideal
subalgebra $\textbf{U}$ that contains all group-like elements has
a PBW-basis $T$ over $\textbf{k}[G]$ which can be extended up to a
PBW-basis $P_T$ of $H$. The set $T$ of PBW-generators for $\textbf{U}$ can be
obtained from the PBW-basis of $H$ given in \cite[Theorem 2]{Kh2} in the
following way.

Suppose that for a given PBW-generator $[u]$ there exists an
element $c \in \textbf{U}$
\begin{equation}\label{geradores}
c=[u]^s+\sum \alpha_iW_i+ \sum \beta_j V_j
\end{equation}
where $W_i$ are the basis words starting with a PBW-generator $[v]<[u]$, $D(W_i)=sD(u)$, and $V_j$ are basis words multiplied from the left by a group-like element with $D(V_j)<sD(u)$. We fix one of the elements with minimal $s$
and denote it by $c_u$. Thus, for every PBW-generator
$[u]$ we have at most one element $c_u$. We define the height
function by the following lemma.

\begin{lem}\label{s1t}
\emph{\cite[Lemma 4.3]{Kh5}} In the representation
\eqref{geradores} of the chosen element $c_u$, either $s=1$ or
$p(u,u)$ is a primitive $t$-th root of $1$ and $s=t$, or (in the
case of positive characteristic) $s=t(char \emph{\textbf{k}})^r$.
\end{lem}

If the height of $[u]$ in $H$ is infinite, then the height of
$c_u$ in $\textbf{U}$ is defined to be infinite as well. If the
height of $[u]$ in $H$ equals $t$ and $p(u,u)$ is a primitive
$t$-th root of $1$, then, due to the above lemma, $s=1$ (note that
in the representation \eqref{geradores} the number $s$ is less
than the height of $[u]$). In this case, the height of $c_u$ in
$\textbf{U}$ is supposed to be $t$ as well. If the characteristic
$l$ is positive and the height of $[u]$ in $H$ equals $tl^r$, then
we define the height of $c_u$ in $\textbf{U}$ to be equal to
$tl^r/s$ (thus, in characteristic zero the height of $c_u$ in
$\textbf{U}$ always equals the height of $[u]$ in $H$).

\begin{pro}\label{pbwU}
\emph{\cite[Proposition 4.4]{Kh5}} An element $c \in H$ belongs to \emph{$\textbf{U}$} if and only if all PBW-generators in the PBW-decomposition of $c$ with respect to $P_T$ belong to $T$. In particular, the set $T$ of all chosen $c_u$
with the above defined height function is a set of
PBW-generators of \emph{$\textbf{U}$} over
\emph{$\textbf{k}[G]$}.
\end{pro}

\begin{defi}\label{U+}
Let $C=\parallel a_{ij} \parallel$ be a generalized Cartan matrix
symmetrizable by $D=diag(d_1,\ldots,d_n)$, $d_ia_{ij}=d_ja_{ji}$.
Denote by $\mathfrak{g}$ a Kac-Moody algebra defined by $C$ (see
\cite{kac}). Suppose that the quantification parameters
$p_{ij}=p(x_i,x_j)=\chi^i(g_j)$ are related by
\begin{equation}\label{pis}
p_{ii}=q^{d_i},\quad p_{ij}p_{ji}=q^{d_ia_{ij}},\quad 1\leq i,j
\leq n.
\end{equation}
The \textit{multiparameter quantization $U^+_q(\mathfrak{g})$} of
the Borel subalgebra $\mathfrak{g^+}$ is a character Hopf algebra generated by $x_1,\ldots,x_n,g_1,\ldots,g_n$ and
defined by Serre relations with the skew brackets (\ref{brac}) in
place of the Lie operation:
\begin{equation}\label{serre}
[[\ldots[[x_i,x_j],x_j],\ldots],x_j]=0,\quad 1\leq i\neq j\leq n,
\end{equation}
where $x_j$ appears $1-a_{ji}$ times.
\end{defi}

We note that relations \eqref{pis} diminish the number of free parameters from $n^2$ to $\frac{n(n-1)}{2}+1$. In particular, if $n=2$ we have just a two-parameter family of quantizations.

By \cite[Theorem 6.1]{Kh1}, the left sides of the relations \eqref{serre} are
skew primitive elements in $G\langle X \rangle$. Therefore, the
ideal generated by these elements is a Hopf ideal, while
$U^+_q(\mathfrak{g})$ indeed has a natural character Hopf algebra
structure.

\begin{defi}
If the multiplicative order $t$ of $q$ is finite, then we define
$u^+_q(\mathfrak{g})$ as $G\langle X\rangle/\Lambda$, where
$\Lambda$ is the biggest Hopf ideal in $G\langle X\rangle^{(2)}$,
which is the set (an ideal) of noncommutative polynomials without
free and linear terms. From \cite[Lemma 2.2]{Kh7}, this is a
$\Gamma^+$-homogeneous ideal. Certainly $\Lambda$ contains all
skew primitive elements of $G\langle X\rangle^{(2)}$ (each one of
them generates a Hopf ideal). Hence, by \cite[Theorem 6.1]{Kh1},
relations (\ref{serre}) are still valid in $u^+_q(\mathfrak{g})$.
\end{defi}

\begin{defi}\label{U-}
Consider a new set of variables
$X^-=\{x_1^-,x_2^-,\ldots,x_n^-\}$. Suppose that an Abelian group
$F$, generated by the elements $f_1,f_2,\ldots,f_n$, acts on the
linear space spanned by $X^-$ so that
$(x_i^-)^{f_j}=p_{ji}^{-1}x_i^-$, where $p_{ij}$ are the same
parameters (\ref{pis}) that define $U^+_q(\mathfrak{g})$. The
relations (\ref{pis}) are invariant under substitutions
$p_{ij}\leftarrow p_{ji}^{-1}$, $q\leftarrow q^{-1}$. This allows
us to define the character Hopf algebra $U^-_q(\mathfrak{g})$ as
$U^+_{q^{-1}}(\mathfrak{g})$ with the characters $\chi^i_-$,
$1\leq i\leq n$ such that $\chi^i_-(f_j)=p_{ji}^{-1}$.
\end{defi}

We may extend the characters $\chi^i$ on $G\times F$ in the
following way $$\chi^i(f_j):=p_{ji}=\chi^j(g_i).$$ Indeed, if
$\prod_kf_k^{m_k}=1$ in $F$, then the application to $x_i^-$
implies $\prod_kp_{ki}^{-m_k}=1$, hence
$\chi^i(\prod_kf_k^{m_k})=\prod_kp_{ki}^{-m_k}$ equals $1$ as
well. In the same way, we may extend the characters $\chi^i_-$ on
$G\times F$ so that
\begin{equation*}\label{char-}
\chi^i_-=(\chi^i)^{-1}
\end{equation*}

In what follows we denote by $H$ the quocient group $(G\times
F)/N$, where $N$ is an arbitrary subgroup with $\chi^i(N)=1$,
$1\leq i\leq n$. For example, if the quantification parameters
satisfy additional symmetry conditions $p_{ij}=p_{ji}$, $1\leq
i,j\leq n$, as in the case for the original Drinfeld-Jimbo and
Lusztig quantifications, then
$\chi^i(g_k^{-1}f_k)=p_{ik}^{-1}p_{ki}=1$, and we may take $N$ to
be the subgroup generated by $g_k^{-1}f_k$, $1\leq k\leq n$. In
this particular case the groups $H$, $G$, $F$ may be identified.

In the general case, without loss of generality we may suppose
that $G,F\subseteq H$. Certainly $\chi^i$, $1\leq i \leq n$ are
characters of $H$ and $H$ still acts on the space spanned by
$X\cup X^-$ by means of these characters and their inverses.

\begin{defi}\label{G2}
Consider the skew group algebra $H\langle X\cup X^- \rangle$ as a
character Hopf algebra:
\begin{equation*}
\Delta(x_i)=x_i\otimes 1+g_i\otimes x_i,\quad
\Delta(x_i^-)=x_i^-\otimes 1+f_i \otimes x_i^-,
\end{equation*}
\begin{equation}\label{gfx}
g^{-1}x_ig=\chi^i(g)\cdot x_i,\quad
g^{-1}x_i^-g=(\chi^i)^{-1}(g)\cdot x_i^-,\quad g\in G,
\end{equation}
\begin{equation}\label{gfx2}
f^{-1}x_if=\chi^i_-(f)\cdot x_i,\quad
f^{-1}x_i^-f=(\chi^i_-)^{-1}(f)\cdot x_i^-,\quad f\in F.
\end{equation}
We define the algebra $U_q(\mathfrak{g})$ as a quotient of
$H\langle X\cup X^- \rangle$ by the following relations:
\begin{equation*}
[[\ldots[[x_i,x_j],x_j],\ldots],x_j]=0,\quad 1\leq i\neq j\leq n,
\end{equation*}
\begin{equation*}
[[\ldots[[x_i^-,x_j^-],x_j^-],\ldots],x_j^-]=0,\quad 1\leq i\neq
j\leq n,
\end{equation*}
where $x_j$ and $x_j^-$ appear $1-a_{ji}$ times and
\begin{equation}\label{serre4}
[x_i,x_j^-]=\delta_i^j(1-g_if_i),\quad 1\leq i,j \leq n
\end{equation}
where the brackets are defined on $H\langle X\cup X^- \rangle$ by
the structure of character Hopf algebra as in (\ref{brac}). Since
due to (\ref{pis}) and \cite[Theorem 6.1]{Kh1}, all polynomials in
the above relations are skew primitive in $H\langle X\cup X^-
\rangle$, they define a Hopf ideal of $H\langle X\cup X^-
\rangle$, that is, the natural homomorphism
$$H\langle X\cup X^- \rangle\rightarrow U_q(\mathfrak{g})$$
defines a Hopf algebra structure on $U_q(\mathfrak{g})$.
\end{defi}

\begin{defi}\label{u}
If $q$ has finite multiplicative order, then $u_q(\mathfrak{g})$
is defined by relations (\ref{serre4}) and $u=0$, $u \in \Lambda$,
$u^-=0$, $u^- \in \Lambda^-$, where $\Lambda, \Lambda^-$ are the
biggest Hopf ideals in $G\langle X\rangle^{(2)}$ and $F\langle
X^-\rangle^{(2)}$, which are sets (ideals) of noncommutative
polynomials without free and linear terms.
\end{defi}

Both algebras $U_q(\mathfrak{g})$ and $u_q(\mathfrak{g})$ are
graded by $\Gamma$, the additive group generated by $\Gamma^+$,
provided that we put $D(x_i^-)=-D(x_i)=-x_i$, $D(H)=0$ since in
this way the relations (\ref{serre4}) become homogeneous.

\begin{cor}\label{homogeneous}
\emph{\cite[Corollary 3.3]{Kh6}} If $q$ is not a root of $1$ and
the Cartan matrix $C=\parallel a_{ij} \parallel$ is of finite
type, then every subalgebra \emph{$\textbf{U}$} of $U_q(\mathfrak{g})$
containing $H$ is $\Gamma$-homogeneous.
\end{cor}

This Corollary is not true for $u_q(\mathfrak{g})$ if $q$ is a root of $1$. In this case, we may have subalgebras that contain $H$ and that are not homogeneous.

From \cite{Kh6}, the subalgebra of $U_q(\mathfrak{g})$ generated by $G$ and
$x_1,x_2,\ldots,x_n$ is isomorphic to $U_q^+(\mathfrak{g})$. Analogously, the subalgebra of
$U_q(\mathfrak{g})$ generated by $F$ and $x_1^-,x_2^-,\ldots,x_n^-$ is isomorphic to
$U_q^-(\mathfrak{g})$. Moreover, we have the following triangular decomposition:
\begin{equation}\label{triaU}
U_q(\mathfrak{g})=U_q^-(\mathfrak{g})\otimes_{\textbf{k}[F]}\textbf{k}[H]\otimes_{\textbf{k}[G]}U_q^+(\mathfrak{g}).
\end{equation}

We also have the same fact for $u_q(\mathfrak{g})$:
\begin{equation}\label{triau}
u_q(\mathfrak{g})=u_q^-(\mathfrak{g})\otimes_{\textbf{k}[F]}\textbf{k}[H]\otimes_{\textbf{k}[G]}u_q^+(\mathfrak{g}).
\end{equation}

Now we have all the necessary definitions and results for the following sections.

\section{The algebras $U_q(\mathfrak{g})$, $u_q(\mathfrak{g})$ when $\mathfrak{g}$ is the simple Lie algebra of type $G_2$}

\quad Now we are going to consider the algebras $U_q(\mathfrak{g})$, when $q$ is not a root of $1$, and $u_q(\mathfrak{g})$, when $q$ has finite multiplicative order $t$, $t>4$, $t \neq 6$, where $\mathfrak{g}$ is the simple Lie algebra of type $G_2$. Let us remember that the algebra $U_q(\mathfrak{g})$ is defined by generators $x_1, x_2, x_1^-, x_2^-,g_1,g_2,f_1,f_2$
and relations
\begin{equation*}\label{rel1}
[[[[x_1,x_2],x_2],x_2],x_2]=0,\quad [x_1,[x_1,x_2]]=0,
\end{equation*}
\begin{equation}\label{rel2}
[[[[x_1^-,x_2^-],x_2^-],x_2^-],x_2^-]=0,\quad
[x_1^-,[x_1^-,x_2^-]]=0,
\end{equation}
\begin{equation*}\label{rel3}
[x_i,x_j^-]=\delta_i^j(1-g_if_i),\quad i,j=1,2,
\end{equation*}
according to Definition \ref{G2}, where $p_{11}=q^3,$ $p_{22}=q$ and $p_{12}p_{21}=q^{-3}$. In this way we have a two-parameter family of quantizations. We also have that the generators $g_1,g_2,f_1,f_2$ are all invertible, commute with each other and
\begin{equation*}
x_ig_j=p_{ij}g_jx_i,\quad x_if_j=p_{ji}f_jx_i,\quad x_i^-g_j=p_{ij}^{-1}g_jx_i^-,\quad x_i^-f_j=p_{ji}^{-1}f_jx_i^-.
\end{equation*}
The algebra $u_q(\mathfrak{g})$ is defined as in Definition \ref{u} and has the same relations \eqref{rel2}.

As we have seen in the preliminaries, the algebras $U_q(\mathfrak{g})$ and $u_q(\mathfrak{g})$ have the triangular decompositions \eqref{triaU} and \eqref{triau}. Now we will give some results that are proved for $U_q^+(\mathfrak{g})$ and $u_q^+(\mathfrak{g})$.

\begin{teo}\label{pbw1} \emph{\cite[Theorem 3.4]{BSP}}
If $q$ is not a root of $1$, then the following elements
\begin{align}\label{base2}
[A]&=x_1,\notag\\ [B]&=[x_1,x_2],\notag\\
[C]&=[[x_1,x_2],[[x_1,x_2],x_2]],\\ [D]&=[[x_1,x_2],x_2],\notag\\
[E]&=[[[x_1,x_2],x_2],x_2],\notag\\ [F]&=x_2.\notag
\end{align}
form a set of PBW-generators for $U_q^+(\mathfrak{g})$ over
\emph{$\textbf{k}[G]$}, and each PBW-generator has infinite height.
If we suppose that $x_1>x_2$, then $A > B > C > D
> E > F$.
\end{teo}

\begin{teo}\label{pbwq}\emph{\cite[Theorem 3.6]{BSP}}
If $q$ has finite multiplicative order $t$, $t>4$, $t\neq 6$, then
the elements from list
\eqref{base2} form a set of PBW-generators for
$u_q^+(\mathfrak{g})$ over \emph{$\textbf{k}[G]$}. The height $h$
of $[u] \in \{[B], [D], [F]\}$ equals $t$. For $[u] \in \{[A],
[C], [E]\}$ we have $h=t$ if $3$ is not a divisor of $t$ and
$h=\frac{t}{3}$ otherwise. In all cases $[u]^h=0$ in
$u_q^+(\mathfrak{g})$.
\end{teo}


\begin{lem}\label{s1}\emph{\cite[Corollary 3.7]{BSP}}
The exponent $s$ given in \eqref{geradores} is $1$ for every $[u]$
in the PBW-basis.
\end{lem}

\begin{teo}\label{lattice} \emph{\cite[Theorem 1.1]{BSP}}
If $q$ is not a root of 1, the lattice of right coideal
subalgebras containing $\textbf{\emph{k}}[G]$ of
$U_q^+(\mathfrak{g})$ is given in the Figure 1. If $q$ has
multiplicative order $t>4$, $t \neq 6$, the same figure is the
lattice of homogeneous right coideal subalgebras containing
$\textbf{\emph{k}}[G]$ of $u_q^+(\mathfrak{g})$.

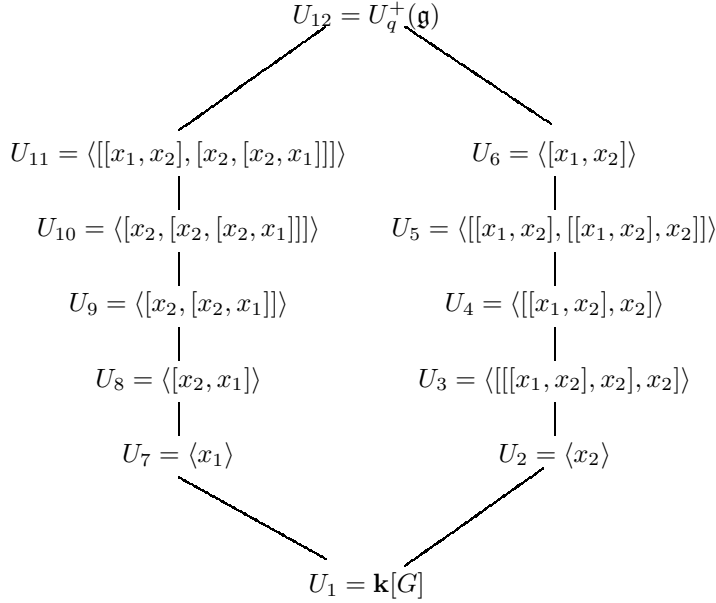
\begin{figure}[!htb]
\begin{center}
\unitlength 1mm
\begin{picture}(60.00,81.25)(0,0)

\put(35.00,81.25){\makebox(0,0)[cc]{$U_{12}=U_q^+(\mathfrak{g})$}}

\put(20.00,71.25){\makebox(0,0)[cc]{}}

\put(60.00,61.25){\makebox(0,0)[bc]{$U_6=\langle[x_1,x_2]\rangle$}}

\put(0.00,5.00){\makebox(0,0)[cc]{}}

\put(60.00,51.25){\makebox(0,0)[bc]{$U_5=\langle[[x_1,x_2],[[x_1,x_2],x_2]]\rangle$}}

\put(60.00,41.25){\makebox(0,0)[bc]{$U_4=\langle[[x_1,x_2],x_2]\rangle$}}

\put(60.00,31.25){\makebox(0,0)[bc]{$U_3=\langle[[[x_1,x_2],x_2],x_2]\rangle$}}

\put(60.00,21.25){\makebox(0,0)[bc]{$U_2=\langle x_{2}\rangle$}}

\put(35.00,5.63){\makebox(0,0)[cc]{$U_1=\textbf{k}[G]$}}

\put(55.00,0.00){\makebox(0,0)[cc]{}}

\linethickness{0.15mm}
\put(60.00,55.63){\line(0,1){4.38}}

\put(10.00,51.25){\makebox(0,0)[bc]{$U_{10}=\langle[x_{2},[x_{2},[x_{2},x_{1}]]]\rangle$}}

\put(10.00,41.25){\makebox(0,0)[bc]{$U_9=\langle[x_{2},[x_{2},x_{1}]]\rangle$}}

\put(10.00,31.25){\makebox(0,0)[bc]{$U_8=\langle[x_{2},x_{1}]\rangle$}}

\put(10.00,21.25){\makebox(0,0)[bc]{$U_7=\langle x_{1}\rangle$}}

\linethickness{0.15mm}


\put(10.00,61.25){\makebox(0,0)[bc]{$U_{11}=\langle[[x_{1},x_{2}],[x_{2},[x_{2},x_{1}]]]\rangle$}}

\linethickness{0.15mm}
\put(60.00,45.63){\line(0,1){4.38}}

\linethickness{0.15mm}
\put(60.00,35.63){\line(0,1){4.38}}

\linethickness{0.15mm}
\put(60.00,25.63){\line(0,1){4.38}}

\linethickness{0.15mm}
\put(10.00,55.63){\line(0,1){4.38}}

\linethickness{0.15mm}
\put(10.00,45.63){\line(0,1){4.38}}

\linethickness{0.15mm}
\put(10.00,35.63){\line(0,1){4.38}}

\linethickness{0.15mm}
\put(10.00,25.63){\line(0,1){4.38}}

\linethickness{0.15mm}
\multiput(40.00,80.00)(0.18,-0.12){109}{\line(1,0){0.18}}

\linethickness{0.15mm}
\multiput(10.00,20.00)(0.18,-0.1){109}{\line(1,0){0.18}}

\linethickness{0.15mm}
\multiput(10.00,66.25)(0.17,0.12){115}{\line(1,0){0.17}}

\linethickness{0.15mm}
\multiput(40.00,8.25)(0.17,0.12){109}{\line(1,0){0.17}}
\end{picture}
\caption{Lattice of Right Coideal Subalgebras}
\end{center}
\end{figure}
\end{teo}

Since the algebra $U_q^-(\mathfrak{g})$ is defined as $U_{q^{-1}}^+(\mathfrak{g})$, a similar result is valid for the negative quantum Borel algebra too.

At this point we would like to say that any ($\Gamma$-homogeneous) right
coideal subalgebra of $U_q(\mathfrak{g})$ (respectively, $u_q(\mathfrak{g})$) that contains
$\textbf{k}[H]$ has the triangular decomposition, and that for any
two right coideal subalgebras
$\textbf{k}[F]\subseteq\textbf{U}^-\subseteq U_q^-(\mathfrak{g})$,
$\textbf{k}[G]\subseteq\textbf{U}^+\subseteq U_q^+(\mathfrak{g})$
(respectively, $\textbf{k}[F]\subseteq\textbf{U}^-\subseteq
u_q^-(\mathfrak{g})$, $\textbf{k}[G]\subseteq\textbf{U}^+\subseteq
u_q^+(\mathfrak{g})$) the tensor product
\begin{equation}\label{triangU}
\textbf{U}=\textbf{U}^-\otimes_{\textbf{k}[F]}\textbf{k}[H]\otimes_{\textbf{k}[G]}\textbf{U}^+
\end{equation}
is a right coideal subalgebra. However this is not true: the tensor product \eqref{triangU} is
not always a subalgebra. We have instead the next two lemmas.

\begin{lem}\label{lema1}
If $q$ is not a root of 1 then every right coideal subalgebra
\emph{$\textbf{U}\supseteq \textbf{k}[H]$} of $U_q(\mathfrak{g})$ has a
decomposition \eqref{triangU}, where \emph{$\textbf{U}^+\supseteq
\textbf{k}[G]$} and \emph{$\textbf{U}^-\supseteq \textbf{k}[F]$}
are right coideal subalgebras of $U_q^+(\mathfrak{g})$ and $U_q^-(\mathfrak{g})$,
respectively. If $q$ has finite multiplicative order $t,$ $t>4,$ $t\neq
6,$, then
this is the case for $\Gamma$-homogeneous right coideal
subalgebras of $u_q(\mathfrak{g})$.
\end{lem}

\proof Due to the triangular decompositions \eqref{triaU} and \eqref{triau}, the set $P=P_-\cup P_+$ is a set of PBW-generators for
$U_q(\mathfrak{g})$ (respectively, $u_q(\mathfrak{g})$) over $\textbf{k}[H]$, where $P_-$ is the
set of PBW-generators for $U_q^-(\mathfrak{g})$ $$\{x_1^-, x_2^-,
[x_1^-,x_2^-], [[x_1^-,x_2^-],x_2^-],
[[[x_1^-,x_2^-],x_2^-],x_2^-],
[[x_1^-,x_2^-],[[x_1^-,x_2^-],x_2^-]]\}$$ and $P_+$ is the set of
PBW-generators for $U_q^+(\mathfrak{g})$ $$\{x_1, x_2, [x_1,x_2],
[[x_1,x_2],x_2], [[[x_1,x_2],x_2],x_2],
[[x_1,x_2],[[x_1,x_2],x_2]]\}.$$

Let us fix the following order on the skew primitive generators
\begin{equation*}
x_1>x_2>x_1^->x_2^-. \end{equation*} By Proposition \ref{pbwU} and
Lemma \ref{s1}, the subalgebra $\textbf{U}$ has PBW-generators of
the form
\begin{equation*} [u]+\sum\alpha_iW_i+\sum\beta_jV_j
\in \textbf{U}, \end{equation*} where $[u] \in P$, $W_i$ are the
basis words starting with PBW-generators smaller than $[u]$,
$D(W_i)=D(u)$, and $V_j$ are basis words multiplied from the left by a group-like element with $D(V_j)<D(u)$.
From the definition of the degree given in \eqref{grau}, all $W_i$
have the same constitution as the leading term $[u]$. Thus all
$W_i$'s and the leading term $[u]$ belong to the same component of
the triangular decomposition (that is, they have only elements
from $X$, or only elements from $X^-$, they can't have both). It
remains to show that there are no terms $V_j$.

If $q$ is not a root of $1$, by Corollary \ref{homogeneous}, the
algebra $\textbf{U}$ is $\Gamma$-homogeneous. Hence (in both
cases) the PBW-generators may be chosen to be $\Gamma$-homogeneous
as well. It means that $a-c=a_j-c_j$ and $b-d=b_j-d_j$ supposing
that $ax_1+bx_2+cx_1^-+dx_2^-$ is the degree of $[u]$ and
$a_jx_1+b_jx_2+c_jx_1^-+d_jx_2^-$ is the degree of $V_j$, for
every $j$. However it contradicts the hypothesis that
$D(V_j)<D(u)$, as it will be shown.

If the leading term $[u] \in P_-$, then $a=b=0$, what provides
$a_j=b_j=0$, from the order defined in \eqref{ordemgrau}. So
$c=c_j$ and $d=d_j$, and $D(V_j)=D(u)$. If $[u] \in P_+$, then
$c=d=0$ and $a=a_j-c_j$, $b=b_j-d_j$. But then $a_j\geq a$,
$b_j\geq b$, $c_j\geq c$ and $d_j\geq d$, what gives us
$D(V_j)\geq D(u)$.

Now we can see that all PBW-generators belong either to
$U_q^-(\mathfrak{g})$ (respectively, $u_q^-(\mathfrak{g})$) or $U_q^+(\mathfrak{g})$ (respectively, $u_q^+(\mathfrak{g})$).
Therefore, $\textbf{U}$ has the decomposition \eqref{triangU}.
\hfill $\square$ \vspace{0.5cm}

\begin{lem}\label{lema2}
\emph{\cite[Lemma 9.3]{Kh6}} The tensor product \eqref{triangU}
is a right coideal subalgebra if and only if
\begin{equation}\label{u+u-}
[\textbf{U}^+,\textbf{U}^-]\subseteq\textbf{U}^-\otimes_{\textbf{k}[F]}\textbf{k}[H]\otimes_{\textbf{k}[G]}\textbf{U}^+.
\end{equation}
\end{lem}

In fact, it was also proved by Kharchenko that, to verify the inclusion \eqref{u+u-}, we just need to verify that
\begin{equation*}
[[u^+],[u^-]]\in\textbf{U}^-\otimes_{\textbf{k}[F]}\textbf{k}[H]\otimes_{\textbf{k}[G]}\textbf{U}^+.
\end{equation*}
for every PBW-generators $[u^+]$ and $[u^-]$ of $\textbf{U}^+$ and $\textbf{U}^-$, respectively.

\section{Right Coideal Subalgebras of $U_q(\mathfrak{g}), u_q(\mathfrak{g})$}

\quad The idea for this section is, using Theorem \ref{lattice}
and Lemma \ref{lema2}, to describe the lattice of (homogeneous) right coideal
subalgebras containing $\textbf{k}[H]$ of $U_q(\mathfrak{g})$ (respectively, $u_q(\mathfrak{g})$). This lattice is represented by four different figures. The first one is the combination of the right coideal subalgebras from the right side of Figure 1 with the right side of the lattice of right coideal subalgebras of $U_q^-(\mathfrak{g})$. The second is the combination of the left sides of both lattices. The last two analyze the left side with the right side, and the right side with the left side. In what follows we denote by $\langle u^-;\, v\rangle$ the space $\langle u^-\rangle\otimes_{\textbf{k}[F]}\textbf{k}[H]\otimes_{\textbf{k}[G]}\langle v\rangle$, where $u,v$ are generators of the right coideal subalgebras of quantum Borel subalgebras from Figure 1. If we get $\langle u_1^-;\, v_1\rangle$, $\langle u_2^-;\, v_2\rangle$ two right coideal subalgebras from two different figures, then their supremum is the smallest right coideal subalgebra containing $\langle u_3^-;\, v_3\rangle$ and their infimum is the greatest right coideal subalgebra contained in $\langle u_0^-;\, v_0\rangle$, where $\langle u^-_3\rangle$ is the supremum of $\langle u^-_1\rangle$ and $\langle u^-_2\rangle$, $\langle v_3\rangle$ is the supremum of $\langle v_1\rangle$ and $\langle v_2\rangle$, $\langle u^-_0\rangle$ is the infimum of $\langle u^-_1\rangle$ and $\langle u^-_2\rangle$ and $\langle v_0\rangle$ is the infimum of $\langle v_1\rangle$ and $\langle v_2\rangle$. Our main result says that the space $\langle u^-;\, v\rangle$ is a right coideal subalgebra if and only if it appears on one of the following four figures.

\begin{teo}\label{teorema}
If $q \in \emph{\textbf{k}}$ is not a root of $1$, the lattice of right coideal subalgebras containing $\textbf{k}[H]$ of $U_q(\mathfrak{g})$ is given by the following four figures. If $q$ has finite multiplicative order $t$, $t>4, t\neq 6$, the same lattice represents the homogeneous right coideal subalgebras containing $\textbf{k}[H]$ of $u_q(\mathfrak{g})$.
\end{teo}





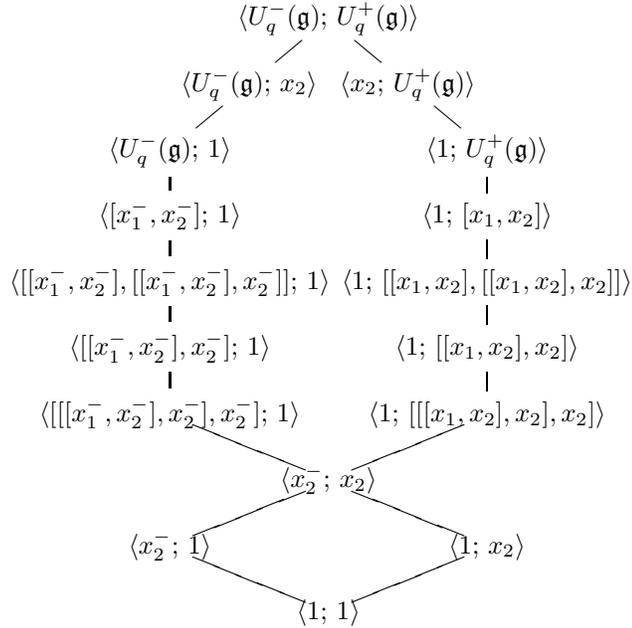
\begin{center}
$\begin{diagram}[tight,width=1.5em,height=1.25em]
  &            &    &            &    \langle U_q^-(\mathfrak{g});\,U_q^+(\mathfrak{g})\rangle      &            &  &             &\\
              &     &           & \ld   &        & \rdLine    &               &     & \\
  &            & \langle U_q^-(\mathfrak{g});\,x_2\rangle     &            &        &            &\langle x_2;\,U_q^+(\mathfrak{g}) \rangle&             &\\
              & \ld    &            &    &        &    &               &\rdLine      & \\
\langle U_q^-(\mathfrak{g});\,1\rangle &            &            &            &        &            &               &             &\langle 1;\,U_q^+(\mathfrak{g})\rangle\\
\dLine &            &            &            &        &            &               &             &\dLine\\
\langle [x_1^-,x_2^-];\,1\rangle &            &            &            &        &            &               &             &\langle 1;\,[x_1,x_2]\rangle\\
\dLine &            &            &            &        &            &               &             &\dLine\\
\langle [[x_1^-,x_2^-],[[x_1^-,x_2^-],x_2^-]];\,1\rangle &            &            &            &        &            &               &             &\langle 1;\,[[x_1,x_2],[[x_1,x_2],x_2]]\rangle\\
\dLine &            &            &            &        &            &               &             &\dLine\\
\langle [[x_1^-,x_2^-],x_2^-];\,1\rangle &            &            &            &        &            &               &             &\langle 1;\,[[x_1,x_2],x_2]\rangle\\
\dLine &            &            &            &        &            &               &             &\dLine\\
\langle [[[x_1^-,x_2^-],x_2^-],x_2^-];\,1\rangle &            &            &            &        &            &               &             &\langle 1;\,[[[x_1,x_2],x_2],x_2]\rangle\\
                     &\rdLine(4,2)&            &            &        &            &               &\ldLine(4,2) &\\
  &            &            &            &\langle x_2^-;\,x_2 \rangle &            &               &             & \\
               &            &            &\ldLine(4,2)&        &\rdLine(4,2)&               &             &\\
\langle x_2^-;\,1\rangle &            &            &            &        &            &               &             &\langle 1;\,x_2\rangle\\
                     &\rdLine(4,2)&            &            &        &            &               &\ldLine(4,2) &\\
                     &            &            &            & \langle 1;\,1 \rangle &            &               &             &
\end{diagram}$
\end{center}

\begin{landscape}
\tiny
$
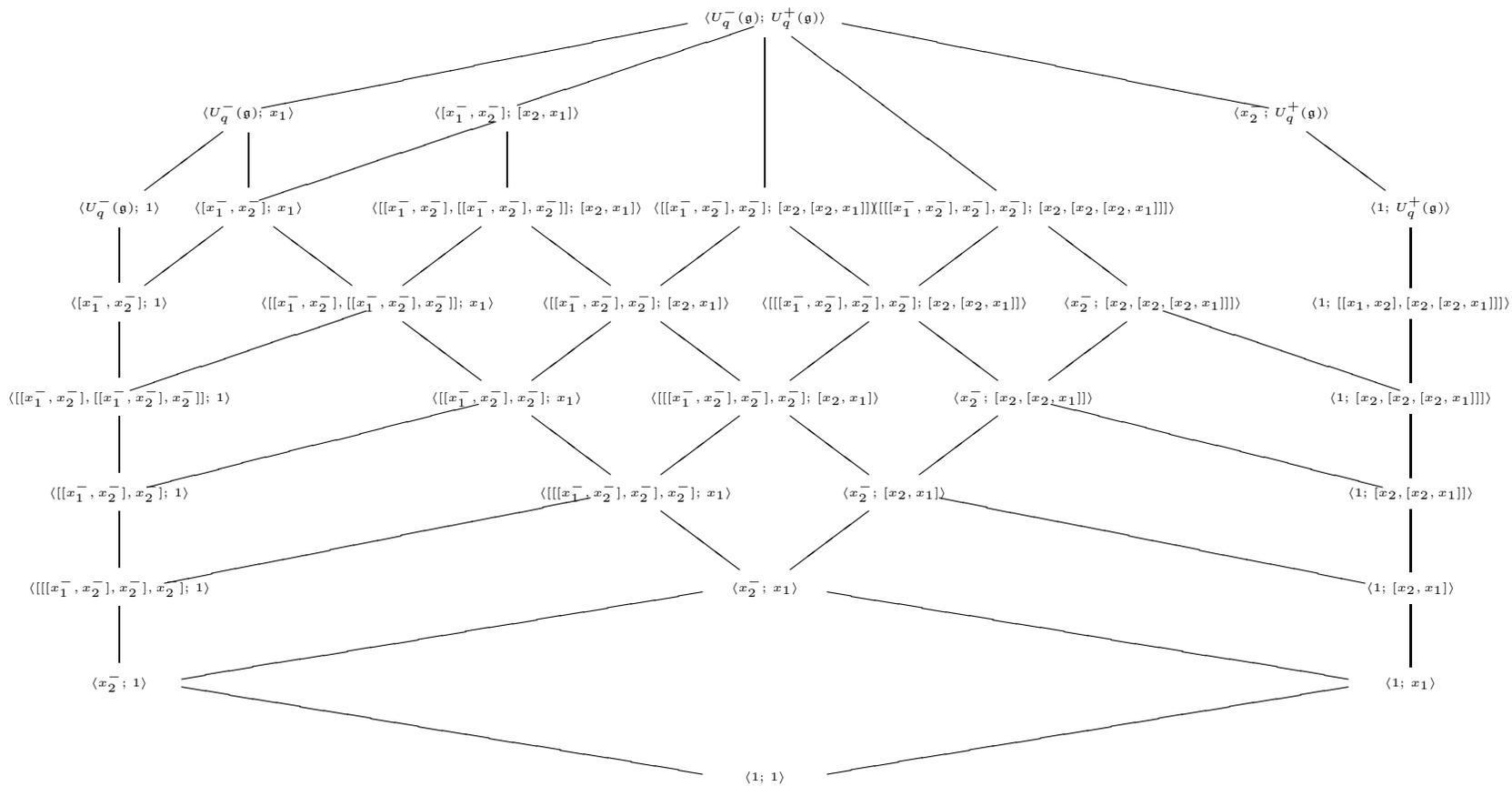
\begin{diagram}[tight,width=3.95em,height=2.9em]
 & & & & & & & & & &\langle U_q^-(\mathfrak{g});\,U_q^+(\mathfrak{g})\rangle & & & & & & & & & & \\
 & & & & & & & & \ld(7,2)&\ld(4,2) &\dd &\rd(4,4) &\rd(7,2) & & & & & & & & \\
 & &\langle U_q^-(\mathfrak{g});\,x_1\rangle & & & &\langle [x_1^-,x_2^-];\,[x_2,x_1]\rangle & & & & & & & & & & & &\langle x_2^-;\,U_q^+(\mathfrak{g})\rangle& & \\
 &\ld & \dd & & &  \ld(4,2) &\dd & & & & & & & & & & & & &\rd & \\
\langle U_q^-(\mathfrak{g});\,1\rangle & &\langle [x_1^-,x_2^-];\,x_1\rangle & & & &\langle [[x_1^-,x_2^-],[[x_1^-,x_2^-],x_2^-]];\,[x_2,x_1]\rangle & & & &\langle [[x_1^-,x_2^-],x_2^-];\,[x_2,[x_2,x_1]]\rangle & & & &\langle [[[x_1^-,x_2^-],x_2^-],x_2^-];\,[x_2,[x_2,[x_2,x_1]]]\rangle & & & & & &\langle 1;\,U_q^+(\mathfrak{g})\rangle \\
\dd &\ld & &\rd & &\ld & &\rd & &\ld & &\rd & &\ld & &\rd & & & & &\dd \\
\langle [x_1^-,x_2^-];\,1\rangle & & & &\langle [[x_1^-,x_2^-],[[x_1^-,x_2^-],x_2^-]];\,x_1\rangle & & & &\langle [[x_1^-,x_2^-],x_2^-];\,[x_2,x_1]\rangle & & & &\langle [[[x_1^-,x_2^-],x_2^-],x_2^-];\,[x_2,[x_2,x_1]]\rangle & & & &\langle x_2^-;\,[x_2,[x_2,[x_2,x_1]]]\rangle & & & &\langle 1;\,[[x_1,x_2],[x_2,[x_2,x_1]]]\rangle \\
\dd &  & &\ld(4,2) & &\rd & &\ld & &\rd & &\ld & &\rd & &\ld & &\rd(4,2) & & &\dd \\
\langle [[x_1^-,x_2^-],[[x_1^-,x_2^-],x_2^-]];\,1\rangle & & & & & &\langle [[x_1^-,x_2^-],x_2^-];\,x_1\rangle & & & &\langle [[[x_1^-,x_2^-],x_2^-],x_2^-];\,[x_2,x_1]\rangle & & & &\langle x_2^-;\,[x_2,[x_2,x_1]]\rangle & & & & & &\langle 1;\,[x_2,[x_2,[x_2,x_1]]]\rangle \\
\dd & & & & &\ld(6,2) & &\rd & &\ld & &\rd & &\ld & &\rd(6,2) & & & & &\dd \\
\langle [[x_1^-,x_2^-],x_2^-];\,1\rangle & & & & & & & &\langle [[[x_1^-,x_2^-],x_2^-],x_2^-];\,x_1\rangle & & & &\langle x_2^-;\,[x_2,x_1]\rangle & & & & & & & & \langle 1;\,[x_2,[x_2,x_1]]\rangle\\
\dd & & & & & & &\ld(8,2) & &\rd & &\ld & &\rd(8,2) & & & & & & &\dd \\
\langle [[[x_1^-,x_2^-],x_2^-],x_2^-];\,1\rangle & & & & & & & & & &\langle x_2^-;\,x_1\rangle & & & & & & & & & & \langle 1;\,[x_2,x_1]\rangle\\
\dd & & & & & & & & &\ld(10,2) & &\rd(10,2) & & & & & & & & &\dd \\
\langle x_2^-;\,1\rangle & & & & & & & & & & & & & & & & & & & & \langle 1;\,x_1\rangle\\
 &\rd(10,2) & & & & & & & & & & & & & & & & & &\ld(10,2) & \\
 & & & & & & & & & &\langle 1;\,1\rangle & & & & & & & & & & \\
\end{diagram}
$

$
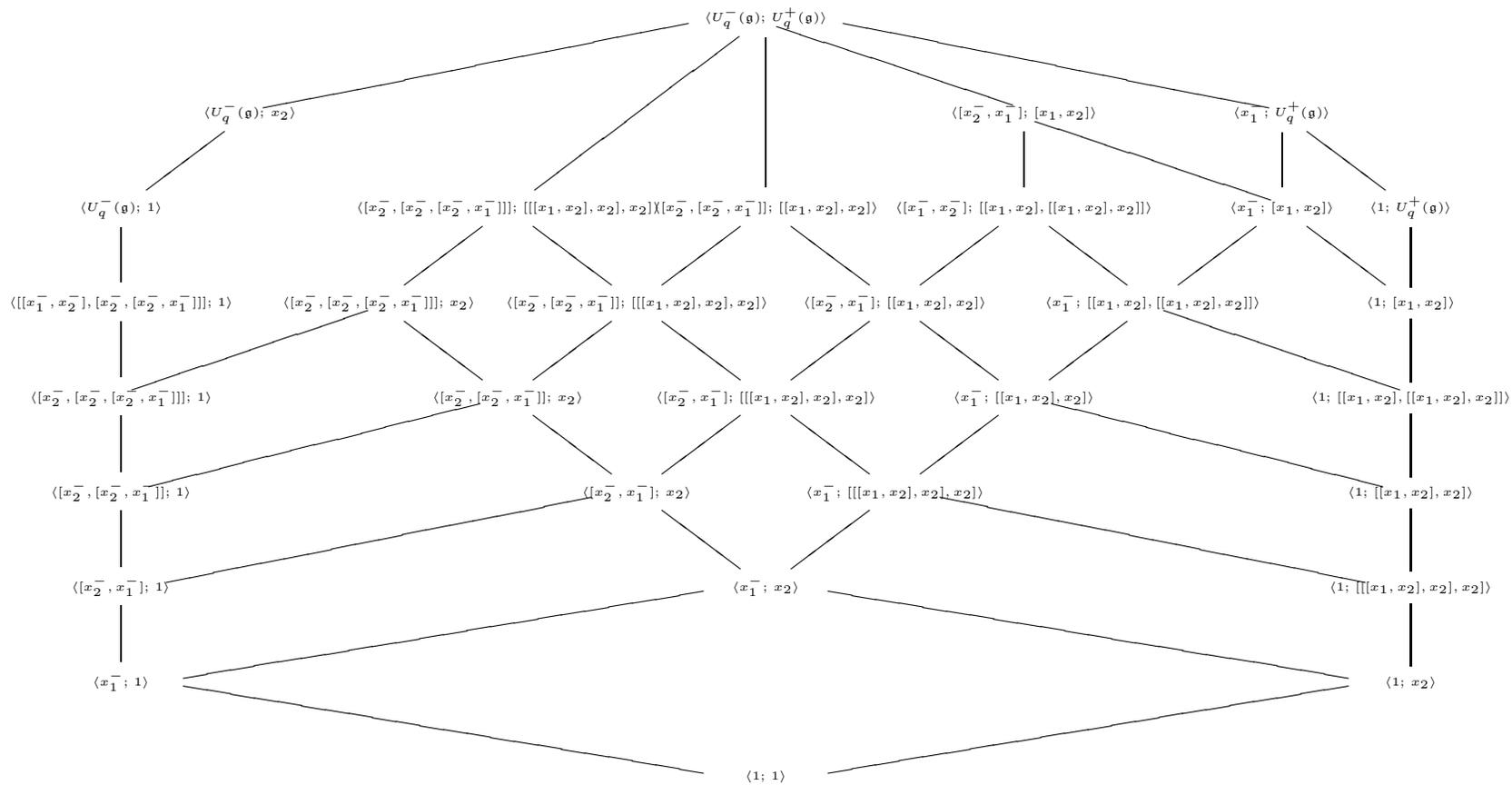
\begin{diagram}[tight,width=3.95em,height=2.9em]
 & & & & & & & & & &\langle U_q^-(\mathfrak{g});\,U_q^+(\mathfrak{g})\rangle & & & & & & & & & & \\
 & & & & & & & & \ld(7,2)&\ld(4,4) &\dd &\rd(4,2) &\rd(7,2) & & & & & & & & \\
 & &\langle U_q^-(\mathfrak{g});\,x_2\rangle & & & & & & & & & & & &\langle [x_2^-,x_1^-];\,[x_1,x_2]\rangle & & & &\langle x_1^-;\,U_q^+(\mathfrak{g})\rangle& & \\
 &\ld & & & & & & & & & & & & &\dd &\rd(4,2) & & &\dd &\rd & \\
\langle U_q^-(\mathfrak{g});\,1\rangle & & & & & &\langle [x_2^-,[x_2^-,[x_2^-,x_1^-]]];\,[[[x_1,x_2],x_2],x_2]\rangle & & & &\langle [x_2^-,[x_2^-,x_1^-]];\,[[x_1,x_2],x_2]\rangle & & & &\langle [x_1^-,x_2^-];\,[[x_1,x_2],[[x_1,x_2],x_2]]\rangle & & & &\langle x_1^-;\,[x_1,x_2]\rangle & &\langle 1;\,U_q^+(\mathfrak{g})\rangle \\
\dd & & & & &\ld & &\rd & &\ld & &\rd & &\ld & &\rd & &\ld & & \rd&\dd \\
\langle [[x_1^-,x_2^-],[x_2^-,[x_2^-,x_1^-]]];\,1\rangle & & & &\langle [x_2^-,[x_2^-,[x_2^-,x_1^-]]];\,x_2\rangle & & & &\langle [x_2^-,[x_2^-,x_1^-]];\,[[[x_1,x_2],x_2],x_2]\rangle & & & &\langle [x_2^-,x_1^-];\,[[x_1,x_2],x_2]\rangle & & & &\langle x_1^-;\,[[x_1,x_2],[[x_1,x_2],x_2]]\rangle & & & &\langle 1;\,[x_1,x_2]\rangle \\
\dd &  & &\ld(4,2) & &\rd & &\ld & &\rd & &\ld & &\rd & &\ld & &\rd(4,2) & & &\dd \\
\langle [x_2^-,[x_2^-,[x_2^-,x_1^-]]];\,1\rangle & & & & & &\langle [x_2^-,[x_2^-,x_1^-]];\,x_2\rangle & & & &\langle [x_2^-,x_1^-];\,[[[x_1,x_2],x_2],x_2]\rangle & & & &\langle x_1^-;\,[[x_1,x_2],x_2]\rangle & & & & & &\langle 1;\,[[x_1,x_2],[[x_1,x_2],x_2]]\rangle \\
\dd & & & & &\ld(6,2) & &\rd & &\ld & &\rd & &\ld & &\rd(6,2) & & & & &\dd \\
\langle [x_2^-,[x_2^-,x_1^-]];\,1\rangle & & & & & & & &\langle [x_2^-,x_1^-];\,x_2\rangle & & & &\langle x_1^-;\,[[[x_1,x_2],x_2],x_2]\rangle & & & & & & & & \langle 1;\,[[x_1,x_2],x_2]\rangle\\
\dd & & & & & & &\ld(8,2) & &\rd & &\ld & &\rd(8,2) & & & & & & &\dd \\
\langle [x_2^-,x_1^-];\,1\rangle & & & & & & & & & &\langle x_1^-;\,x_2\rangle & & & & & & & & & & \langle 1;\,[[[x_1,x_2],x_2],x_2]\rangle\\
\dd & & & & & & & & &\ld(10,2) & &\rd(10,2) & & & & & & & & &\dd \\
\langle x_1^-;\,1\rangle & & & & & & & & & & & & & & & & & & & & \langle 1;\,x_2\rangle\\
 &\rd(10,2) & & & & & & & & & & & & & & & & & &\ld(10,2) & \\
 & & & & & & & & & &\langle 1;\,1\rangle & & & & & & & & & & \\
\end{diagram}
$
\end{landscape}

\normalsize
\begin{center}
$
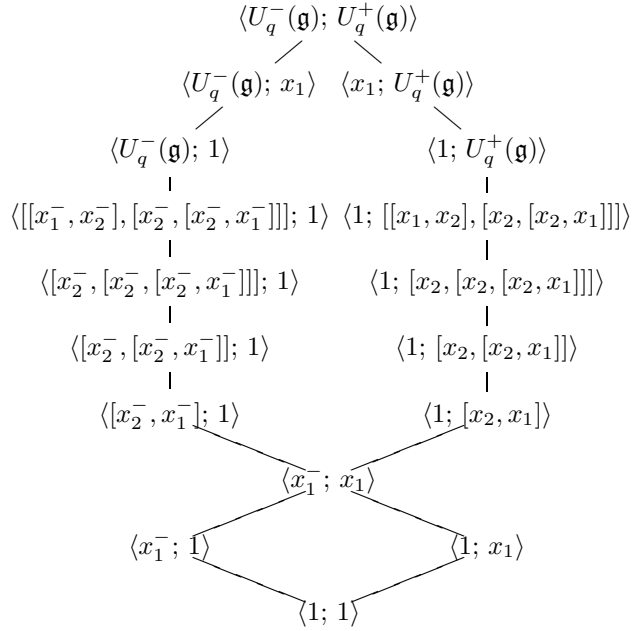
\begin{diagram}[tight,width=1.5em,height=1.25em]
  &            &    &            &    \langle U_q^-(\mathfrak{g});\,U_q^+(\mathfrak{g})\rangle      &            &  &             &\\
              &     &           & \ld   &        & \rdLine    &               &     & \\
  &            & \langle U_q^-(\mathfrak{g});\,x_1\rangle     &            &        &            &\langle x_1;\,U_q^+(\mathfrak{g}) \rangle&             &\\
              & \ld    &            &    &        &    &               &\rdLine      & \\
\langle U_q^-(\mathfrak{g});\,1\rangle &            &            &            &        &            &               &             &\langle 1;\,U_q^+(\mathfrak{g})\rangle\\
\dLine &            &            &            &        &            &               &             &\dLine\\
\langle [[x_1^-,x_2^-],[x_2^-,[x_2^-,x_1^-]]];\,1\rangle &            &            &            &        &            &               &             &\langle 1;\,[[x_1,x_2],[x_2,[x_2,x_1]]]\rangle\\
\dLine &            &            &            &        &            &               &             &\dLine\\
\langle [x_2^-,[x_2^-,[x_2^-,x_1^-]]];\,1\rangle &            &            &            &        &            &               &             &\langle 1;\,[x_2,[x_2,[x_2,x_1]]]\rangle\\
\dLine &            &            &            &        &            &               &             &\dLine\\
\langle [x_2^-,[x_2^-,x_1^-]];\,1\rangle &            &            &            &        &            &               &             &\langle 1;\,[x_2,[x_2,x_1]]\rangle\\
\dLine &            &            &            &        &            &               &             &\dLine\\
\langle [x_2^-,x_1^-];\,1\rangle &            &            &            &        &            &               &             &\langle 1;\,[x_2,x_1]\rangle\\
                     &\rdLine(4,2)&            &            &        &            &               &\ldLine(4,2) &\\
  &            &            &            &\langle x_1^-;\,x_1 \rangle &            &               &             & \\
               &            &            &\ldLine(4,2)&        &\rdLine(4,2)&               &             &\\
\langle x_1^-;\,1\rangle &            &            &            &        &            &               &             &\langle 1;\,x_1\rangle\\
                     &\rdLine(4,2)&            &            &        &            &               &\ldLine(4,2) &\\
                     &            &            &            & \langle 1;\,1 \rangle &            &               &             &
\end{diagram}$
\end{center}

\proof
From Lemma \ref{lema1}, all (homogeneous) right coideal subalgebras containing $\textbf{k}[H]$ of $U_q(\mathfrak{g})$ (respectively, $u_q(\mathfrak{g})$) have the form
\begin{equation}
\textbf{U}=\textbf{U}^-\otimes_{\textbf{k}[F]}\textbf{k}[H]\otimes_{\textbf{k}[G]}\textbf{U}^+,
\end{equation}
where $\textbf{U}^+,\textbf{U}^-$ are (homogeneous) right coideal subalgebras containing $\textbf{k}[G]$,$\textbf{k}[F]$ of $U_q^+(\mathfrak{g})$,$U_q^-(\mathfrak{g})$ (respectively, $u_q^+(\mathfrak{g})$,$u_q^-(\mathfrak{g})$). Now we just have to see which combinations are right coideal subalgebras and which ones are not.

First let us remember that, from Theorem \ref{lattice} the (homogeneous)
right coideal subalgebras of $U_q^+(\mathfrak{g})$ (respectively, $u_q^+(\mathfrak{g})$) with respective
PBW-generators over $\textbf{k}[G]$ are:

\begin{itemize}

\item $U_1=\textbf{k}[G]$

\item $U_2=\langle x_2\rangle$

$\{x_2\}$

\item $U_3=\langle [[[x_1,x_2],x_2],x_2]\rangle$

$\{x_2;\, [[[x_1,x_2],x_2],x_2]\}$

\item $U_4=\langle [[x_1,x_2],x_2]\rangle$

$\{x_2;\, [[[x_1,x_2],x_2],x_2];\, [[x_1,x_2],x_2]\}$

\item $U_5=\langle [[x_1,x_2],[[x_1,x_2],x_2]]\rangle$

$\{x_2;\, [[[x_1,x_2],x_2],x_2];\, [[x_1,x_2],x_2];\,
[[x_1,x_2],[[x_1,x_2],x_2]]\}$

\item $U_6=\langle [x_1,x_2]\rangle$

$\{x_2;\, [[[x_1,x_2],x_2],x_2];\, [[x_1,x_2],x_2];\,
[[x_1,x_2],[[x_1,x_2],x_2]];\, [x_1,x_2]\}$

\item $U_7=\langle x_1\rangle$

$\{x_1\}$

\item $U_8=\langle [x_2,x_1]\rangle$

$\{x_1;\, [x_2,x_1]\}$

\item $U_9=\langle [x_2,[x_2,x_1]]\rangle$

$\{x_1;\, [x_2,x_1];\, [x_2,[x_2,x_1]]\}$

\item $U_{10}=\langle [x_2,[x_2,[x_2,x_1]]]\rangle$

$\{x_1;\, [x_2,x_1];\, [x_2,[x_2,x_1]];\, [x_2,[x_2,[x_2,x_1]]]\}$

\item $U_{11}=\langle [[x_1,x_2],[x_2,[x_2,x_1]]]\rangle$

$\{x_1;\, [x_2,x_1];\, [x_2,[x_2,x_1]];\, [x_2,[x_2,[x_2,x_1]]];\,
[[x_1,x_2],[x_2,[x_2,x_1]]]\}$

\item $U_{12}=U_q^+(\mathfrak{g})$

$\{x_2;\, [[[x_1,x_2],x_2],x_2];\, [[x_1,x_2],x_2];\,
[[x_1,x_2],[[x_1,x_2],x_2]];\, [x_1,x_2];\, x_1\}$

\end{itemize}

Since $U_q^-(\mathfrak{g})=U_{q^{-1}}^+(\mathfrak{g})$, the right coideal subalgebras containing $\textbf{k}[F]$ of $U_q^-(\mathfrak{g})$ are exactly the
same, just substituting $x_1$ for $x_1^-$ and $x_2$ for $x_2^-$.
From Lemmas \ref{lema1} and \ref{lema2}, in order to find out the
(homogeneous) right coideal subalgebras containing $\textbf{k}[H]$ of $U_q(\mathfrak{g})$ (respectively, $u_q(\mathfrak{g})$) we just have to calculate
the skew-commutators between the possible PBW-generators, which are just listed here. Since we do not need to know the exact coefficients, for simplification of the expressions we define the following relation
 $$u \sim v \quad \text{if and only if} \quad u=\alpha v,$$
where $\alpha \in \textbf{k}, \alpha \neq 0$. For the coefficients $\alpha_i, \,i=1,\ldots,8$, that appear on some expressions, we also have $\alpha_i \in \textbf{k}$ and $\alpha_i \neq 0$, except for $\alpha_1$ which may be zero. All calculations will be shown at the appendix, where these coefficients are also specified. \label{13}

\begin{description}
\item [(A1)] $[x_1,x_1^-]\sim 1-g_1f_1$,

\item [(A2)] $[x_1,x_2^-]\sim 0$,

\item [(A3)]
$[x_1,[x_1^-,x_2^-]]\sim x_2^-g_1f_1$,

\item [(A4)]
$[x_1,[x_2^-,x_1^-]]\sim x_2^-$,

\item [(A5)]
$[x_1,[[x_1^-,x_2^-],x_2^-]]\sim (x_2^-)^2g_1f_1$,

\item [(A6)]
$[x_1,[[[x_1^-,x_2^-],x_2^-],x_2^-]]\sim (x_2^-)^3g_1f_1$,

\item [(A7)]
$[x_1,[[x_1^-,x_2^-],[[x_1^-,x_2^-],x_2^-]]]\sim x_2^-[[x_1^-,x_2^-],x_2^-]g_1f_1+\alpha_1[[[x_1^-,x_2^-],x_2^-],x_2^-]g_1f_1$,

\item [(A8)] $[x_2,x_1^-]\sim 0$,

\item [(A9)] $[x_2,x_2^-]\sim 1-g_2f_2$,

\item [(A10)] $[x_2,[x_1^-,x_2^-]]\sim x_1^-$,

\item [(A11)] $[x_2,[x_2^-,x_1^-]]\sim x_1^-g_2f_2$,

\item [(A12)]
$[x_2,[[x_1^-,x_2^-],x_2^-]]\sim [x_1^-,x_2^-]$,

\item [(A13)]
$[x_2,[x_2^-,[x_2^-,x_1^-]]]\sim [x_2^-,x_1^-]g_2f_2$,

\item [(A14)]
$[x_2,[[[x_1^-,x_2^-],x_2^-],x_2^-]]\sim [[x_1^-,x_2^-],x_2^-]$,

\item [(A15)]
$[x_2,[x_2^-,[x_2^-,[x_2^-,x_1^-]]]]\sim [x_2^-,[x_2^-,x_1^-]]g_2f_2$,

\item [(A16)]
$[x_2,[[x_1^-,x_2^-],[x_2^-,[x_2^-,x_1^-]]]]\sim [x_2^-,x_1^-]^2+\alpha_2[x_2^-,x_1^-]^2g_2f_2+\alpha_3[x_2^-,[x_2^-,x_1^-]]x_1^-g_2f_2+\alpha_4x_2^-[x_2^-,x_1^-]x_1^-g_2f_2$,

\item [(A17)]
$[[x_1,x_2],x_1^-]\sim g_1f_1x_2$,

\item [(A18)]
$[[x_1,x_2],x_2^-]\sim x_1$,

\item [(A19)]
$[[x_1,x_2],[x_2^-,x_1^-]]\sim 1-g_1f_1g_2f_2$,

\item [(A20)]
$[[x_1,x_2],[x_2^-,[x_2^-,x_1^-]]]\sim x_2^-$,

\item [(A21)] $[[x_2,x_1],x_1^-]\sim x_2$,

\item [(A22)] $[[x_2,x_1],x_2^-]\sim g_2f_2x_1$,

\item [(A23)] $[[x_2,x_1],[x_1^-,x_2^-]]\sim 1-g_2f_2g_1f_1$,

\item [(A24)] $[[x_2,x_1],[[x_1^-,x_2^-],x_2^-]]\sim x_2^-g_2f_2g_1f_1$,

\item [(A25)] $[[x_2,x_1],[[[x_1^-,x_2^-],x_2^-],x_2^-]]\sim (x_2^-)^2g_2f_2g_1f_1$,

\item [(A26)] $[[x_2,x_1],[[x_1^-,x_2^-],[[x_1^-,x_2^-],x_2^-]]]\sim [[x_1^-,x_2^-],x_2^-]g_2f_2g_1f_1$,

\item [(A27)]
$[[[x_1,x_2],x_2],x_1^-]\sim g_1f_1x_2^2$,

\item [(A28)]
$[[[x_1,x_2],x_2],x_2^-]\sim [x_1,x_2]$,

\item [(A29)]
$[[[x_1,x_2],x_2],[x_2^-,x_1^-]]\sim g_1f_1g_2f_2x_2$,

\item [(A30)]
$[[[x_1,x_2],x_2],[x_2^-,[x_2^-,x_1^-]]]\sim 1-g_1f_1g_2^2f_2^2$,

\item [(A31)]
$[[[x_1,x_2],x_2],[x_2^-,[x_2^-,[x_2^-,x_1^-]]]]\sim x_2^-$,

\item [(A32)]
$[[x_2,[x_2,x_1]],x_1^-]\sim x_2^2$,

\item [(A33)]
$[[x_2,[x_2,x_1]],x_2^-]\sim g_2f_2[x_2,x_1]$,

\item [(A34)]
$[[x_2,[x_2,x_1]],[x_1^-,x_2^-]]\sim x_2$,

\item [(A35)]
$[[x_2,[x_2,x_1]],[[x_1^-,x_2^-],x_2^-]]\sim 1-g_1f_1g_2^2f_2^2$,

\item [(A36)]
$[[x_2,[x_2,x_1]],[[[x_1^-,x_2^-],x_2^-],x_2^-]]\sim x_2^-g_1f_1g_2^2f_2^2$,

\item[(A37)]
$[[x_2,[x_2,x_1]],[[x_1^-,x_2^-],[[x_1^-,x_2^-],x_2^-]]]\sim [x_1^-,x_2^-]$,

\item [(A38)]
$[[[[x_1,x_2],x_2],x_2],x_1^-]\sim g_1f_1x_2^3$,

\item [(A39)]
$[[[[x_1,x_2],x_2],x_2],x_2^-]\sim [[x_1,x_2],x_2]$,

\item [(A40)]
$[[[[x_1,x_2],x_2],x_2],[x_2^-,x_1^-]]\sim g_1f_1g_2f_2x_2^2$,

\item [(A41)]
$[[[[x_1,x_2],x_2],x_2],[x_2^-,[x_2^-,x_1^-]]]\sim g_1f_1g_2^2f_2^2x_2$,

\item [(A42)]
$[[[[x_1,x_2],x_2],x_2],[x_2^-,[x_2^-,[x_2^-,x_1^-]]]]\sim (1-g_1f_1g_2f_2)g_2f_2$,

\item [(A43)]
$[[x_2,[x_2,[x_2,x_1]]],x_2^-]\sim g_2f_2[x_2,[x_2,x_1]]$,

\item [(A44)]
$[[x_2,[x_2,[x_2,x_1]]],[[x_1^-,x_2^-],x_2^-]]\sim x_2$,

\item [(A45)]
$[[x_2,[x_2,[x_2,x_1]]],[[[x_1^-,x_2^-],x_2^-],x_2^-]]\sim 1-g_1f_1g_2^3f_2^3$,

\item [(A46)]
$[[[x_1,x_2],[x_2,[x_2,x_1]]],x_2^-]\sim [x_2,x_1]^2+\alpha_5g_2f_2[x_2,x_1]^2+\alpha_6g_2f_2x_2[x_2,x_1]x_1+\alpha_7g_2f_2[x_2,[x_2,x_1]]x_1$,

\item [(A47)]
$[[[x_1,x_2],[[x_1,x_2],x_2]],x_2^-]\sim [x_1,x_2]^2$,

\item [(A48)]
$[[[x_1,x_2],[[x_1,x_2],x_2]],x_1^-]\sim g_1f_1[[[x_1,x_2],x_2],x_2]+\alpha_8g_1f_1x_2[[x_1,x_2],x_2]$,

\item [(A49)]
$[[[x_1,x_2],[[x_1,x_2],x_2]],[x_2^-,x_1^-]]\sim g_1f_1g_2f_2[[x_1,x_2],x_2]$,

\item [(A50)]
$[[[x_1,x_2],[[x_1,x_2],x_2]],[x_2^-,[x_2^-,x_1^-]]]\sim [x_1,x_2]$.

\end{description}

Now we are ready to apply Lemma \ref{lema2} and see which spaces $\langle u^-;\, v\rangle$ are right coideal subalgebras and which ones are not.

 Since $U_q^+(\mathfrak{g})$ and $U_q^-(\mathfrak{g})$ (respectively, $u_q^+(\mathfrak{g})$ and $u_q^-(\mathfrak{g})$) have 12 (homogeneous) right coideal subalgebras containing $\textbf{k}[G]$ and $\textbf{k}[F]$, from Lemma \ref{lema1} we have 144 possible (homogeneous) right coideal subalgebras to study.

 From the formula $ug=\chi^u(g)gu$, $g \in G$, we obtain that $$[\textbf{k}[G],\textbf{U}^-]=\textbf{k}[G]\textbf{U}^--\textbf{U}^-\textbf{k}[G]=\textbf{U}^-\textbf{k}[G].$$ Analogously,  $[\textbf{U}^+,\textbf{k}[F]]=\textbf{k}[F]\textbf{U}^+$, and we have 23 right coideal subalgebras of the form $\textbf{k}[F]\otimes_{\textbf{k}[F]}\textbf{k}[H]\otimes_{\textbf{k}[G]}\textbf{U}^{+}$ or $\textbf{U}^{-}\otimes_{\textbf{k}[F]}\textbf{k}[H]\otimes_{\textbf{k}[G]}\textbf{k}[G]$, where $\textbf{U}^{\pm}$ runs through the set of all right coideal subalgebras of quantum Borel subalgebras.

From \textbf{(A1)}, page \pageref{13}, we get $$[x_1,x_1^-]=1-g_1f_1\in \langle x_1^- \rangle\otimes_{\textbf{k}[F]}\textbf{k}[H]\otimes_{\textbf{k}[G]}\langle x_1 \rangle.$$

Using \textbf{(A8)} we have
\begin{align*}
[x_2,x_1^-]=0 &\in \langle x_1^- \rangle\otimes_{\textbf{k}[F]}\textbf{k}[H]\otimes_{\textbf{k}[G]}\langle x_2 \rangle\\&\subseteq \langle x_1^- \rangle\otimes_{\textbf{k}[F]}\textbf{k}[H]\otimes_{\textbf{k}[G]}\langle [[[x_1,x_2],x_2],x_2] \rangle\\&\subseteq \langle x_1^- \rangle\otimes_{\textbf{k}[F]}\textbf{k}[H]\otimes_{\textbf{k}[G]}\langle [[x_1,x_2],x_2] \rangle\\&\subseteq \langle x_1^- \rangle\otimes_{\textbf{k}[F]}\textbf{k}[H]\otimes_{\textbf{k}[G]}\langle [[x_1,x_2],[[x_1,x_2],x_2]] \rangle\\&\subseteq \langle x_1^- \rangle\otimes_{\textbf{k}[F]}\textbf{k}[H]\otimes_{\textbf{k}[G]}\langle [x_1,x_2] \rangle\\&\subseteq \langle x_1^- \rangle\otimes_{\textbf{k}[F]}\textbf{k}[H]\otimes_{\textbf{k}[G]}U_q^+(\mathfrak{g}).
\end{align*}

From \textbf{(A38)},
\begin{align*}
[[[[x_1,x_2],x_2],x_2],x_1^-]&\sim g_1f_1x_2^3\\&\in \langle x_1^- \rangle\otimes_{\textbf{k}[F]}\textbf{k}[H]\otimes_{\textbf{k}[G]}\langle [[[x_1,x_2],x_2],x_2] \rangle\\&\subseteq \langle x_1^- \rangle\otimes_{\textbf{k}[F]}\textbf{k}[H]\otimes_{\textbf{k}[G]}\langle [[x_1,x_2],x_2] \rangle\\&\subseteq \langle x_1^- \rangle\otimes_{\textbf{k}[F]}\textbf{k}[H]\otimes_{\textbf{k}[G]}\langle [[x_1,x_2],[[x_1,x_2],x_2]] \rangle\\&\subseteq \langle x_1^- \rangle\otimes_{\textbf{k}[F]}\textbf{k}[H]\otimes_{\textbf{k}[G]}\langle [x_1,x_2] \rangle\\&\subseteq \langle x_1^- \rangle\otimes_{\textbf{k}[F]}\textbf{k}[H]\otimes_{\textbf{k}[G]}U_q^+(\mathfrak{g}).
\end{align*}

With \textbf{(A27)},
\begin{align*}
[[[x_1,x_2],x_2],x_1^-]&\sim g_1f_1x_2^2\\&\in \langle x_1^- \rangle\otimes_{\textbf{k}[F]}\textbf{k}[H]\otimes_{\textbf{k}[G]}\langle [[x_1,x_2],x_2] \rangle\\&\subseteq \langle x_1^- \rangle\otimes_{\textbf{k}[F]}\textbf{k}[H]\otimes_{\textbf{k}[G]}\langle [[x_1,x_2],[[x_1,x_2],x_2]] \rangle\\&\subseteq \langle x_1^- \rangle\otimes_{\textbf{k}[F]}\textbf{k}[H]\otimes_{\textbf{k}[G]}\langle [x_1,x_2] \rangle\\&\subseteq \langle x_1^- \rangle\otimes_{\textbf{k}[F]}\textbf{k}[H]\otimes_{\textbf{k}[G]}U_q^+(\mathfrak{g}).
\end{align*}

Using \textbf{(A48)}, \begin{align*}
[[[x_1,x_2],[[x_1,x_2],x_2]],x_1^-]&\sim g_1f_1[[[x_1,x_2],x_2],x_2]+\alpha g_1f_1x_2[[x_1,x_2],x_2]\\&\in \langle x_1^- \rangle\otimes_{\textbf{k}[F]}\textbf{k}[H]\otimes_{\textbf{k}[G]}\langle [[x_1,x_2],[[x_1,x_2],x_2]] \rangle\\&\subseteq \langle x_1^- \rangle\otimes_{\textbf{k}[F]}\textbf{k}[H]\otimes_{\textbf{k}[G]}\langle [x_1,x_2] \rangle\\&\subseteq \langle x_1^- \rangle\otimes_{\textbf{k}[F]}\textbf{k}[H]\otimes_{\textbf{k}[G]}U_q^+(\mathfrak{g}).\end{align*}

The equation \textbf{(A17)} provides
\begin{align*}[[x_1,x_2],x_1^-]\sim g_1f_1x_2&\in \langle x_1^- \rangle\otimes_{\textbf{k}[F]}\textbf{k}[H]\otimes_{\textbf{k}[G]}\langle [x_1,x_2] \rangle\\&\subseteq \langle x_1^- \rangle\otimes_{\textbf{k}[F]}\textbf{k}[H]\otimes_{\textbf{k}[G]}U_q^+(\mathfrak{g}).\end{align*}
Then, we have 7 right coideal subalgebras:
\begin{description}
  \item [Y24.] $\langle x_1^- \rangle\otimes_{\textbf{k}[F]}\textbf{k}[H]\otimes_{\textbf{k}[G]}\langle x_1 \rangle$,
  \item [Y25.] $\langle x_1^- \rangle\otimes_{\textbf{k}[F]}\textbf{k}[H]\otimes_{\textbf{k}[G]}\langle x_2 \rangle$,
  \item [Y26.] $\langle x_1^- \rangle\otimes_{\textbf{k}[F]}\textbf{k}[H]\otimes_{\textbf{k}[G]}\langle [[[x_1,x_2],x_2],x_2] \rangle$,
  \item [Y27.] $\langle x_1^- \rangle\otimes_{\textbf{k}[F]}\textbf{k}[H]\otimes_{\textbf{k}[G]}\langle [[x_1,x_2],x_2] \rangle$,
  \item [Y28.] $\langle x_1^- \rangle\otimes_{\textbf{k}[F]}\textbf{k}[H]\otimes_{\textbf{k}[G]}\langle [[x_1,x_2],[[x_1,x_2],x_2]] \rangle$,
  \item [Y29.] $\langle x_1^- \rangle\otimes_{\textbf{k}[F]}\textbf{k}[H]\otimes_{\textbf{k}[G]}\langle [x_1,x_2] \rangle$,
  \item [Y30.] $\langle x_1^- \rangle\otimes_{\textbf{k}[F]}\textbf{k}[H]\otimes_{\textbf{k}[G]}U_q^+(\mathfrak{g})$.
\end{description}

Analogously, from \textbf{(A2)}, \textbf{(A9)}, \textbf{(A18)}, \textbf{(A22)}, \textbf{(A28)}, \textbf{(A33)}, \textbf{(A39)}, \textbf{(A43)} and \textbf{(A47)} we have another 6 right coideal subalgebras:
\begin{description}
  \item [Y31.] $\langle x_2^- \rangle\otimes_{\textbf{k}[F]}\textbf{k}[H]\otimes_{\textbf{k}[G]}\langle x_2 \rangle$,
  \item [Y32.] $\langle x_2^- \rangle\otimes_{\textbf{k}[F]}\textbf{k}[H]\otimes_{\textbf{k}[G]}\langle x_1 \rangle$,
  \item [Y33.] $\langle x_2^- \rangle\otimes_{\textbf{k}[F]}\textbf{k}[H]\otimes_{\textbf{k}[G]}\langle [x_2,x_1] \rangle$,
  \item [Y34.] $\langle x_2^- \rangle\otimes_{\textbf{k}[F]}\textbf{k}[H]\otimes_{\textbf{k}[G]}\langle [x_2,[x_2,x_1]] \rangle$,
  \item [Y35.] $\langle x_2^- \rangle\otimes_{\textbf{k}[F]}\textbf{k}[H]\otimes_{\textbf{k}[G]}\langle [x_2,[x_2,[x_2,x_1]]] \rangle$,
  \item [Y36.] $\langle x_2^- \rangle\otimes_{\textbf{k}[F]}\textbf{k}[H]\otimes_{\textbf{k}[G]}U_q^+(\mathfrak{g})$.
\end{description}

Using \textbf{(A11)}, \textbf{(A19)}, \textbf{(A29)}, \textbf{(A40)} and \textbf{(A49)}, and the already analyzed cases, we obtain 5 more right coideal subalgebras:

\begin{description}
  \item [Y37.] $\langle [x_2^-,x_1^-] \rangle\otimes_{\textbf{k}[F]}\textbf{k}[H]\otimes_{\textbf{k}[G]}\langle x_2 \rangle$,
  \item [Y38.] $\langle [x_2^-,x_1^-] \rangle\otimes_{\textbf{k}[F]}\textbf{k}[H]\otimes_{\textbf{k}[G]}\langle [[[x_1,x_2],x_2],x_2] \rangle$,
  \item [Y39.] $\langle [x_2^-,x_1^-] \rangle\otimes_{\textbf{k}[F]}\textbf{k}[H]\otimes_{\textbf{k}[G]}\langle [[x_1,x_2],x_2] \rangle$,
  \item [Y40.] $\langle [x_2^-,x_1^-] \rangle\otimes_{\textbf{k}[F]}\textbf{k}[H]\otimes_{\textbf{k}[G]}\langle [[x_1,x_2],[[x_1,x_2],x_2]] \rangle$,
  \item [Y41.] $\langle [x_2^-,x_1^-] \rangle\otimes_{\textbf{k}[F]}\textbf{k}[H]\otimes_{\textbf{k}[G]}\langle [x_1,x_2] \rangle$.
\end{description}

With \textbf{(A13)}, \textbf{(A30)} and \textbf{(A41)} we also have:

\begin{description}
  \item [Y42.] $\langle [x_2^-,[x_2^-,x_1^-]] \rangle\otimes_{\textbf{k}[F]}\textbf{k}[H]\otimes_{\textbf{k}[G]}\langle x_2 \rangle$,
  \item [Y43.] $\langle [x_2^-,[x_2^-,x_1^-]] \rangle\otimes_{\textbf{k}[F]}\textbf{k}[H]\otimes_{\textbf{k}[G]}\langle [[[x_1,x_2],x_2],x_2] \rangle$,
  \item [Y44.] $\langle [x_2^-,[x_2^-,x_1^-]] \rangle\otimes_{\textbf{k}[F]}\textbf{k}[H]\otimes_{\textbf{k}[G]}\langle [[x_1,x_2],x_2] \rangle$.
\end{description}

From \textbf{(A15)} and \textbf{(A42)} we get 2 new right coideal subalgebras:

\begin{description}
  \item [Y45.] $\langle [x_2^-,[x_2^-,[x_2^-,x_1^-]]] \rangle\otimes_{\textbf{k}[F]}\textbf{k}[H]\otimes_{\textbf{k}[G]}\langle x_2 \rangle$,
  \item [Y46.] $\langle [x_2^-,[x_2^-,[x_2^-,x_1^-]]] \rangle\otimes_{\textbf{k}[F]}\textbf{k}[H]\otimes_{\textbf{k}[G]}\langle [[[x_1,x_2],x_2],x_2] \rangle$.
\end{description}

Using \textbf{(A6)}, \textbf{(A25)}, \textbf{(A36)} and \textbf{(A45)} we obtain the following right coideal subalgebras:

\begin{description}
  \item [Y47.] $\langle [[[x_1^-,x_2^-],x_2^-],x_2^-] \rangle\otimes_{\textbf{k}[F]}\textbf{k}[H]\otimes_{\textbf{k}[G]}\langle x_1 \rangle$,
  \item [Y48.] $\langle [[[x_1^-,x_2^-],x_2^-],x_2^-] \rangle\otimes_{\textbf{k}[F]}\textbf{k}[H]\otimes_{\textbf{k}[G]}\langle [x_2,x_1] \rangle$,
  \item [Y49.] $\langle [[[x_1^-,x_2^-],x_2^-],x_2^-] \rangle\otimes_{\textbf{k}[F]}\textbf{k}[H]\otimes_{\textbf{k}[G]}\langle [x_2,[x_2,x_1]] \rangle$,
  \item [Y50.] $\langle [[[x_1^-,x_2^-],x_2^-],x_2^-] \rangle\otimes_{\textbf{k}[F]}\textbf{k}[H]\otimes_{\textbf{k}[G]}\langle [x_2,[x_2,[x_2,x_1]]] \rangle$.
\end{description}

With \textbf{(A5)}, \textbf{(A24)} and \textbf{(A35)} we have also:
\begin{description}
  \item [Y51.] $\langle [[x_1^-,x_2^-],x_2^-] \rangle \otimes_{\textbf{k}[F]}\textbf{k}[H]\otimes_{\textbf{k}[G]}\langle x_1 \rangle,$
  \item [Y52.] $\langle [[x_1^-,x_2^-],x_2^-] \rangle \otimes_{\textbf{k}[F]}\textbf{k}[H]\otimes_{\textbf{k}[G]}\langle [x_2,x_1] \rangle$,
  \item [Y53.] $\langle [[x_1^-,x_2^-],x_2^-] \rangle \otimes_{\textbf{k}[F]}\textbf{k}[H]\otimes_{\textbf{k}[G]}\langle [x_2,[x_2,x_1]] \rangle$.
\end{description}

From \textbf{(A7)} and \textbf{(A26)} we get that the following possibilities are actually right coideal subalgebras:

\begin{description}
  \item [Y54.] $\langle [[x_1^-,x_2^-],[[x_1^-,x_2^-],x_2^-]] \rangle\otimes_{\textbf{k}[F]}\textbf{k}[H]\otimes_{\textbf{k}[G]}\langle x_1 \rangle$,
  \item [Y55.] $\langle [[x_1^-,x_2^-],[[x_1^-,x_2^-],x_2^-]] \rangle\otimes_{\textbf{k}[F]}\textbf{k}[H]\otimes_{\textbf{k}[G]}\langle [x_2,x_1] \rangle$.
\end{description}

Using \textbf{(A3)} and \textbf{(A23)} we obtain:

\begin{description}
\item [Y56.] $\langle [x_1^-,x_2^-] \rangle\otimes_{\textbf{k}[F]}\textbf{k}[H]\otimes_{\textbf{k}[G]}\langle x_1 \rangle$,

\item [Y57.] $\langle [x_1^-,x_2^-] \rangle\otimes_{\textbf{k}[F]}\textbf{k}[H]\otimes_{\textbf{k}[G]}\langle [x_2,x_1] \rangle$.
\end{description}

Finally, from \textbf{(A1)} to \textbf{(A3)}, \textbf{(A5)} to \textbf{(A10)}, \textbf{(A12)}, \textbf{(A14)} and \textbf{(A16)}, we have:
\begin{description}
\item [Y58.] $U_q^-(\mathfrak{g})\otimes_{\textbf{k}[F]}\textbf{k}[H]\otimes_{\textbf{k}[G]}\langle x_1 \rangle,$

\item [Y59.] $U_q^-(\mathfrak{g})\otimes_{\textbf{k}[F]}\textbf{k}[H]\otimes_{\textbf{k}[G]}\langle x_2 \rangle.$
\end{description}

We have 59 right coideal subalgebras so far. Adding the trivial case $U_q(\mathfrak{g})=U_q^-(\mathfrak{g})\otimes_{\textbf{k}[F]}\textbf{k}[H]\otimes_{\textbf{k}[G]}U_q^+(\mathfrak{g})$ we complete 60. With these right coideal subalgebras we obtain the 4 suggested lattices. Notice that some right coideal subalgebras may appear on two or more of those four figures. For example, $U_q^-(\mathfrak{g})\otimes_{\textbf{k}[F]}\textbf{k}[H]\otimes_{\textbf{k}[G]}U_q^+(\mathfrak{g})$ is in all four figures, and $U_q^-(\mathfrak{g})\otimes_{\textbf{k}[F]}\textbf{k}[H]\otimes_{\textbf{k}[G]}\langle x_1 \rangle$ is in two.

To complete the proof of this theorem, we still have 84 cases to consider. Now we will see that these cases are not right coideal subalgebras.

From \textbf{(A21)}, $[[x_2,x_1],x_1^-]\sim x_2$ and
\begin{align*}
x_2 \notin &\langle [[x_1,x_2],[x_2,[x_2,x_1]]] \rangle \supseteq\\ &\langle [x_2,[x_2,[x_2,x_1]]] \rangle \supseteq\\ &\langle [x_2,[x_2,x_1]] \rangle \supseteq \langle [x_2,x_1] \rangle,
\end{align*}
so we exclude the following 24 cases:
\begin{description}
\item [N1.] $\langle x_1^-\rangle \otimes_{\textbf{k}[F]}\textbf{k}[H]\otimes_{\textbf{k}[G]}\langle [x_2,x_1] \rangle$,
\item [N2.] $\langle x_1^-\rangle \otimes_{\textbf{k}[F]}\textbf{k}[H]\otimes_{\textbf{k}[G]}\langle [x_2,[x_2,x_1]] \rangle$,
\item [N3.] $\langle x_1^-\rangle \otimes_{\textbf{k}[F]}\textbf{k}[H]\otimes_{\textbf{k}[G]}\langle [x_2,[x_2,[x_2,x_1]]] \rangle$,
\item [N4.] $\langle x_1^-\rangle \otimes_{\textbf{k}[F]}\textbf{k}[H]\otimes_{\textbf{k}[G]}\langle [[x_1,x_2],[x_2,[x_2,x_1]]] \rangle$,

\item [N5.] $\langle [x_2^-,x_1^-]\rangle \otimes_{\textbf{k}[F]}\textbf{k}[H]\otimes_{\textbf{k}[G]}\langle [x_2,x_1] \rangle$,
\item [N6.] $\langle [x_2^-,x_1^-]\rangle \otimes_{\textbf{k}[F]}\textbf{k}[H]\otimes_{\textbf{k}[G]}\langle [x_2,[x_2,x_1]] \rangle$,
\item [N7.] $\langle [x_2^-,x_1^-]\rangle \otimes_{\textbf{k}[F]}\textbf{k}[H]\otimes_{\textbf{k}[G]}\langle [x_2,[x_2,[x_2,x_1]]] \rangle$,
\item [N8.] $\langle [x_2^-,x_1^-]\rangle \otimes_{\textbf{k}[F]}\textbf{k}[H]\otimes_{\textbf{k}[G]}\langle [[x_1,x_2],[x_2,[x_2,x_1]]] \rangle$,

\item [N9.] $\langle [x_2^-,[x_2^-,x_1^-]]\rangle \otimes_{\textbf{k}[F]}\textbf{k}[H]\otimes_{\textbf{k}[G]}\langle [x_2,x_1] \rangle$,
\item [N10.] $\langle [x_2^-,[x_2^-,x_1^-]]\rangle \otimes_{\textbf{k}[F]}\textbf{k}[H]\otimes_{\textbf{k}[G]}\langle [x_2,[x_2,x_1]] \rangle$,
\item [N11.] $\langle [x_2^-,[x_2^-,x_1^-]]\rangle \otimes_{\textbf{k}[F]}\textbf{k}[H]\otimes_{\textbf{k}[G]}\langle [x_2,[x_2,[x_2,x_1]]] \rangle$,
\item [N12.] $\langle [x_2^-,[x_2^-,x_1^-]]\rangle \otimes_{\textbf{k}[F]}\textbf{k}[H]\otimes_{\textbf{k}[G]}\langle [[x_1,x_2],[x_2,[x_2,x_1]]] \rangle$,

\item [N13.] $\langle [x_2^-,[x_2^-,[x_2^-,x_1^-]]]\rangle \otimes_{\textbf{k}[F]}\textbf{k}[H]\otimes_{\textbf{k}[G]}\langle [x_2,x_1] \rangle$,
\item [N14.] $\langle [x_2^-,[x_2^-,[x_2^-,x_1^-]]]\rangle \otimes_{\textbf{k}[F]}\textbf{k}[H]\otimes_{\textbf{k}[G]}\langle [x_2,[x_2,x_1]] \rangle$,
\item [N15.] $\langle [x_2^-,[x_2^-,[x_2^-,x_1^-]]]\rangle \otimes_{\textbf{k}[F]}\textbf{k}[H]\otimes_{\textbf{k}[G]}\langle [x_2,[x_2,[x_2,x_1]]] \rangle$,
\item [N16.] $\langle [x_2^-,[x_2^-,[x_2^-,x_1^-]]]\rangle \otimes_{\textbf{k}[F]}\textbf{k}[H]\otimes_{\textbf{k}[G]}\langle [[x_1,x_2],[x_2,[x_2,x_1]]] \rangle$,

\item [N17.] $\langle [[x_1^-,x_2^-],[x_2^-,[x_2^-,x_1^-]]]\rangle \otimes_{\textbf{k}[F]}\textbf{k}[H]\otimes_{\textbf{k}[G]}\langle [x_2,x_1] \rangle$,
\item [N18.] $\langle [[x_1^-,x_2^-],[x_2^-,[x_2^-,x_1^-]]]\rangle \otimes_{\textbf{k}[F]}\textbf{k}[H]\otimes_{\textbf{k}[G]}\langle [x_2,[x_2,x_1]] \rangle$,
\item [N19.] $\langle [[x_1^-,x_2^-],[x_2^-,[x_2^-,x_1^-]]]\rangle \otimes_{\textbf{k}[F]}\textbf{k}[H]\otimes_{\textbf{k}[G]}\langle [x_2,[x_2,[x_2,x_1]]] \rangle$,
\item [N20.] $\langle [[x_1^-,x_2^-],[x_2^-,[x_2^-,x_1^-]]]\rangle \otimes_{\textbf{k}[F]}\textbf{k}[H]\otimes_{\textbf{k}[G]}\langle [[x_1,x_2],[x_2,[x_2,x_1]]] \rangle$,

\item [N21.] $U_q^-(\mathfrak{g}) \otimes_{\textbf{k}[F]}\textbf{k}[H]\otimes_{\textbf{k}[G]}\langle [x_2,x_1] \rangle$,
\item [N22.] $U_q^-(\mathfrak{g}) \otimes_{\textbf{k}[F]}\textbf{k}[H]\otimes_{\textbf{k}[G]}\langle [x_2,[x_2,x_1]] \rangle$,
\item [N23.] $U_q^-(\mathfrak{g}) \otimes_{\textbf{k}[F]}\textbf{k}[H]\otimes_{\textbf{k}[G]}\langle [x_2,[x_2,[x_2,x_1]]] \rangle$,
\item [N24.] $U_q^-(\mathfrak{g}) \otimes_{\textbf{k}[F]}\textbf{k}[H]\otimes_{\textbf{k}[G]}\langle [[x_1,x_2],[x_2,[x_2,x_1]]] \rangle$.

\end{description}

From \textbf{(A39)}, $[[[[x_1,x_2],x_2],x_2],x_2^-]\sim [[x_1,x_2],x_2]$. So, since we have $$[[x_1,x_2],x_2] \notin \langle [[[[x_1,x_2],x_2],x_2] \rangle,$$ we exclude two more cases:
\begin{description}
\item [N25.] $\langle x_2^- \rangle \otimes_{\textbf{k}[F]}\textbf{k}[H]\otimes_{\textbf{k}[G]}\langle [[[[x_1,x_2],x_2],x_2] \rangle,$
\item [N26.] $U_q^-(\mathfrak{g}) \otimes_{\textbf{k}[F]}\textbf{k}[H]\otimes_{\textbf{k}[G]}\langle [[[[x_1,x_2],x_2],x_2] \rangle.$
\end{description}

Using \textbf{(A28)}, we obtain $[[[x_1,x_2],x_2],x_2^-]\sim [x_1,x_2]$. From the fact that $$[x_1,x_2] \notin \langle [[x_1,x_2],[[x_1,x_2],x_2]]\rangle \supseteq \langle [[x_1,x_2],x_2]\rangle,$$ we eliminate the following possibilities:
\begin{description}
\item [N27.] $\langle x_2^- \rangle \otimes_{\textbf{k}[F]}\textbf{k}[H]\otimes_{\textbf{k}[G]}\langle [[[x_1,x_2],x_2] \rangle,$
\item [N28.] $\langle x_2^- \rangle \otimes_{\textbf{k}[F]}\textbf{k}[H]\otimes_{\textbf{k}[G]}\langle [[[x_1,x_2],x_2] \rangle,$
\item [N29.] $U_q^-(\mathfrak{g}) \otimes_{\textbf{k}[F]}\textbf{k}[H]\otimes_{\textbf{k}[G]}\langle [[x_1,x_2],[[x_1,x_2],x_2]] \rangle,$
\item [N30.] $U_q^-(\mathfrak{g}) \otimes_{\textbf{k}[F]}\textbf{k}[H]\otimes_{\textbf{k}[G]}\langle [[x_1,x_2],[[x_1,x_2],x_2]] \rangle.$
\end{description}

With \textbf{(A18)} we have
$[[x_1,x_2],x_2^-]\sim x_1$, excluding:
\begin{description}
\item [N31.] $\langle x_2^- \rangle \otimes_{\textbf{k}[F]}\textbf{k}[H]\otimes_{\textbf{k}[G]}\langle [x_1,x_2] \rangle,$
\item [N32.] $U_q^-(\mathfrak{g}) \otimes_{\textbf{k}[F]}\textbf{k}[H]\otimes_{\textbf{k}[G]}\langle [x_1,x_2] \rangle.$
\end{description}

From \textbf{(A14)},
$[x_2,[[[x_1^-,x_2^-],x_2^-],x_2^-]]\sim [[x_1^-,x_2^-],x_2^-]$, what eliminates:
\begin{description}
\item [N33.] $\langle [[[[x_1^-,x_2^-],x_2^-],x_2^-] \rangle \otimes_{\textbf{k}[F]}\textbf{k}[H]\otimes_{\textbf{k}[G]}\langle x_2 \rangle,$
\item [N34.] $\langle [[[[x_1^-,x_2^-],x_2^-],x_2^-] \rangle \otimes_{\textbf{k}[F]}\textbf{k}[H]\otimes_{\textbf{k}[G]}\langle [[[x_1,x_2],x_2],x_2] \rangle,$
\item [N35.] $\langle [[[[x_1^-,x_2^-],x_2^-],x_2^-] \rangle \otimes_{\textbf{k}[F]}\textbf{k}[H]\otimes_{\textbf{k}[G]}\langle [[x_1,x_2],x_2] \rangle,$
\item [N36.] $\langle [[[[x_1^-,x_2^-],x_2^-],x_2^-] \rangle \otimes_{\textbf{k}[F]}\textbf{k}[H]\otimes_{\textbf{k}[G]}\langle [[x_1,x_2],[[x_1,x_2],x_2]] \rangle,$
\item [N37.] $\langle [[[[x_1^-,x_2^-],x_2^-],x_2^-] \rangle \otimes_{\textbf{k}[F]}\textbf{k}[H]\otimes_{\textbf{k}[G]}\langle [x_1,x_2] \rangle,$
\item [N38.] $\langle [[[[x_1^-,x_2^-],x_2^-],x_2^-] \rangle \otimes_{\textbf{k}[F]}\textbf{k}[H]\otimes_{\textbf{k}[G]}U_q^+(\mathfrak{g}).$
\end{description}

In \textbf{(A12)} we have $[x_2,[[x_1^-,x_2^-],x_2^-]]\sim [x_1^-,x_2^-]$, and we exclude:
\begin{description}
\item [N39.] $\langle [[[x_1^-,x_2^-],x_2^-] \rangle \otimes_{\textbf{k}[F]}\textbf{k}[H]\otimes_{\textbf{k}[G]}\langle x_2 \rangle,$
\item [N40.] $\langle [[[x_1^-,x_2^-],x_2^-] \rangle \otimes_{\textbf{k}[F]}\textbf{k}[H]\otimes_{\textbf{k}[G]}\langle [[[x_1,x_2],x_2],x_2] \rangle,$
\item [N41.] $\langle [[[x_1^-,x_2^-],x_2^-] \rangle \otimes_{\textbf{k}[F]}\textbf{k}[H]\otimes_{\textbf{k}[G]}\langle [[x_1,x_2],x_2] \rangle,$
\item [N42.] $\langle [[[x_1^-,x_2^-],x_2^-] \rangle \otimes_{\textbf{k}[F]}\textbf{k}[H]\otimes_{\textbf{k}[G]}\langle [[x_1,x_2],[[x_1,x_2],x_2]] \rangle,$
\item [N43.] $\langle [[[x_1^-,x_2^-],x_2^-] \rangle \otimes_{\textbf{k}[F]}\textbf{k}[H]\otimes_{\textbf{k}[G]}\langle [x_1,x_2] \rangle,$
\item [N44.] $\langle [[[x_1^-,x_2^-],x_2^-] \rangle \otimes_{\textbf{k}[F]}\textbf{k}[H]\otimes_{\textbf{k}[G]}U_q^+(\mathfrak{g}),$
\item [N45.] $\langle [[x_1^-,x_2^-],[[[x_1^-,x_2^-],x_2^-]] \rangle \otimes_{\textbf{k}[F]}\textbf{k}[H]\otimes_{\textbf{k}[G]}\langle x_2 \rangle,$
\item [N46.] $\langle [[x_1^-,x_2^-],[[[x_1^-,x_2^-],x_2^-]] \rangle \otimes_{\textbf{k}[F]}\textbf{k}[H]\otimes_{\textbf{k}[G]}\langle [[[x_1,x_2],x_2],x_2] \rangle,$
\item [N47.] $\langle [[x_1^-,x_2^-],[[[x_1^-,x_2^-],x_2^-]] \rangle \otimes_{\textbf{k}[F]}\textbf{k}[H]\otimes_{\textbf{k}[G]}\langle [[x_1,x_2],x_2] \rangle,$
\item [N48.] $\langle [[x_1^-,x_2^-],[[[x_1^-,x_2^-],x_2^-]] \rangle \otimes_{\textbf{k}[F]}\textbf{k}[H]\otimes_{\textbf{k}[G]}\langle [[x_1,x_2],[[x_1,x_2],x_2]] \rangle,$
\item [N49.] $\langle [[x_1^-,x_2^-],[[[x_1^-,x_2^-],x_2^-]] \rangle \otimes_{\textbf{k}[F]}\textbf{k}[H]\otimes_{\textbf{k}[G]}\langle [x_1,x_2] \rangle,$
\item [N50.] $\langle [[x_1^-,x_2^-],[[[x_1^-,x_2^-],x_2^-]] \rangle \otimes_{\textbf{k}[F]}\textbf{k}[H]\otimes_{\textbf{k}[G]}U_q^+(\mathfrak{g}).$
\end{description}

Using \textbf{(A10)}, we know that $[x_2,[x_1^-,x_2^-]]\sim x_1^-$ and eliminate 6 more possibilities:
\begin{description}
\item [N51.] $\langle [x_1^-,x_2^-] \rangle \otimes_{\textbf{k}[F]}\textbf{k}[H]\otimes_{\textbf{k}[G]}\langle x_2 \rangle,$
\item [N52.] $\langle [x_1^-,x_2^-] \rangle \otimes_{\textbf{k}[F]}\textbf{k}[H]\otimes_{\textbf{k}[G]}\langle [[[x_1,x_2],x_2],x_2] \rangle,$
\item [N53.] $\langle [x_1^-,x_2^-] \rangle \otimes_{\textbf{k}[F]}\textbf{k}[H]\otimes_{\textbf{k}[G]}\langle [[x_1,x_2],x_2] \rangle,$
\item [N54.] $\langle [x_1^-,x_2^-] \rangle \otimes_{\textbf{k}[F]}\textbf{k}[H]\otimes_{\textbf{k}[G]}\langle [[x_1,x_2],[[x_1,x_2],x_2]] \rangle,$
\item [N55.] $\langle [x_1^-,x_2^-] \rangle \otimes_{\textbf{k}[F]}\textbf{k}[H]\otimes_{\textbf{k}[G]}\langle [x_1,x_2] \rangle,$
\item [N56.] $\langle [x_1^-,x_2^-] \rangle \otimes_{\textbf{k}[F]}\textbf{k}[H]\otimes_{\textbf{k}[G]}U_q^+(\mathfrak{g}).$
\end{description}

The equality given by \textbf{(A4)} implies
$[x_1,[x_2^-,x_1^-]]\sim x_2^-$, excluding the cases:
\begin{description}
\item [N57.] $\langle [x_2^-,x_1^-] \rangle \otimes_{\textbf{k}[F]}\textbf{k}[H]\otimes_{\textbf{k}[G]}\langle x_1 \rangle,$
\item [N58.] $\langle [x_2^-,x_1^-] \rangle \otimes_{\textbf{k}[F]}\textbf{k}[H]\otimes_{\textbf{k}[G]}U_q^+(\mathfrak{g}),$
\item [N59.] $\langle [x_2^-,[x_2^-,x_1^-]] \rangle \otimes_{\textbf{k}[F]}\textbf{k}[H]\otimes_{\textbf{k}[G]}\langle x_1 \rangle,$
\item [N60.] $\langle [x_2^-,[x_2^-,x_1^-]] \rangle \otimes_{\textbf{k}[F]}\textbf{k}[H]\otimes_{\textbf{k}[G]}U_q^+(\mathfrak{g}),$
\item [N61.] $\langle [x_2^-,[x_2^-,[x_2^-,x_1^-]]] \rangle \otimes_{\textbf{k}[F]}\textbf{k}[H]\otimes_{\textbf{k}[G]}\langle x_1 \rangle,$
\item [N62.] $\langle [x_2^-,[x_2^-,[x_2^-,x_1^-]]] \rangle \otimes_{\textbf{k}[F]}\textbf{k}[H]\otimes_{\textbf{k}[G]}U_q^+(\mathfrak{g}),$
\item [N63.] $\langle [[x_1^-,x_2^-],[x_2^-,[x_2^-,x_1^-]]] \rangle \otimes_{\textbf{k}[F]}\textbf{k}[H]\otimes_{\textbf{k}[G]}\langle x_1 \rangle,$
\item [N64.] $\langle [[x_1^-,x_2^-],[x_2^-,[x_2^-,x_1^-]]] \rangle \otimes_{\textbf{k}[F]}\textbf{k}[H]\otimes_{\textbf{k}[G]}U_q^+(\mathfrak{g}).$
\end{description}

From \textbf{(A44)} we know that
$[[x_2,[x_2,[x_2,x_1]]],[[x_1^-,x_2^-],x_2^-]]\sim x_2$, what eliminates:
\begin{description}
\item [N65.] $\langle [[x_1^-,x_2^-],x_2^-] \rangle \otimes_{\textbf{k}[F]}\textbf{k}[H]\otimes_{\textbf{k}[G]}\langle [x_2,[x_2,[x_2,x_1]]] \rangle,$
\item [N66.] $\langle [[x_1^-,x_2^-],x_2^-] \rangle \otimes_{\textbf{k}[F]}\textbf{k}[H]\otimes_{\textbf{k}[G]}\langle [[x_1,x_2],[x_2,[x_2,x_1]]]\rangle.$
\end{description}

With \textbf{(A37)} we have $[[x_2,[x_2,x_1]],[[x_1^-,x_2^-],[[x_1^-,x_2^-],x_2^-]]]\sim[x_1^-,x_2^-]$, and we obtain that
\begin{description}
\item [N67.] $\langle [[x_1^-,x_2^-],[[x_1^-,x_2^-],x_2^-]] \rangle \otimes_{\textbf{k}[F]}\textbf{k}[H]\otimes_{\textbf{k}[G]}\langle [x_2,[x_2,x_1]] \rangle,$
\item [N68.] $\langle [[x_1^-,x_2^-],[[x_1^-,x_2^-],x_2^-]] \rangle \otimes_{\textbf{k}[F]}\textbf{k}[H]\otimes_{\textbf{k}[G]}\langle [x_2,[x_2,[x_2,x_1]]] \rangle,$
\item [N69.] $\langle [[x_1^-,x_2^-],[[x_1^-,x_2^-],x_2^-]] \rangle \otimes_{\textbf{k}[F]}\textbf{k}[H]\otimes_{\textbf{k}[G]}\langle [[x_1,x_2],[x_2,[x_2,x_1]]]\rangle,$
\end{description}
are not right coideal subalgebras.

From the equality in \textbf{(A34)}, $[[x_2,[x_2,x_1]],[x_1^-,x_2^-]]\sim x_2$, and the following cases are excluded:
\begin{description}
\item [N70.] $\langle [x_1^-,x_2^-] \rangle \otimes_{\textbf{k}[F]}\textbf{k}[H]\otimes_{\textbf{k}[G]}\langle [x_2,[x_2,x_1]] \rangle,$
\item [N71.] $\langle [x_1^-,x_2^-] \rangle \otimes_{\textbf{k}[F]}\textbf{k}[H]\otimes_{\textbf{k}[G]}\langle [x_2,[x_2,[x_2,x_1]]] \rangle,$
\item [N72.] $\langle [x_1^-,x_2^-] \rangle \otimes_{\textbf{k}[F]}\textbf{k}[H]\otimes_{\textbf{k}[G]}\langle [[x_1,x_2],[x_2,[x_2,x_1]]]\rangle.$
\end{description}

In \textbf{(A31)} we have $[[[x_1,x_2],x_2],[x_2^-,[x_2^-,[x_2^-,x_1^-]]]]\sim x_2^-$, eliminating the cases:
\begin{description}
\item [N73.] $\langle [x_2^-,[x_2^-,[x_2^-,x_1^-]]] \rangle \otimes_{\textbf{k}[F]}\textbf{k}[H]\otimes_{\textbf{k}[G]}\langle [[x_1,x_2],x_2]\rangle,$
\item [N74.] $\langle [[x_1^-,x_2^-],[x_2^-,[x_2^-,x_1^-]]] \rangle \otimes_{\textbf{k}[F]}\textbf{k}[H]\otimes_{\textbf{k}[G]}\langle [[x_1,x_2],x_2]\rangle.$
\end{description}

Using \textbf{(A50)} we obtain $[[[x_1,x_2],[[x_1,x_2],x_2]],[x_2^-,[x_2^-,x_1^-]]]\sim [x_1,x_2]$, what excludes:
\begin{description}
\item [N75.] $\langle [x_2^-,[x_2^-,x_1^-]] \rangle \otimes_{\textbf{k}[F]}\textbf{k}[H]\otimes_{\textbf{k}[G]}\langle [[x_1,x_2],[[x_1,x_2],x_2]]\rangle,$
\item [N76.] $\langle [x_2^-,[x_2^-,[x_2^-,x_1^-]]] \rangle \otimes_{\textbf{k}[F]}\textbf{k}[H]\otimes_{\textbf{k}[G]}\langle [[x_1,x_2],[[x_1,x_2],x_2]]\rangle,$
\item [N77.] $\langle [[x_1^-,x_2^-],[x_2^-,[x_2^-,x_1^-]]] \rangle \otimes_{\textbf{k}[F]}\textbf{k}[H]\otimes_{\textbf{k}[G]}\langle [[x_1,x_2],[[x_1,x_2],x_2]]\rangle.$
\end{description}

With \textbf{(A20)}, we have $[[x_1,x_2],[x_2^-,[x_2^-,x_1^-]]]\sim x_2^-$, and we eliminate:
\begin{description}
\item [N78.] $\langle [x_2^-,[x_2^-,x_1^-]] \rangle \otimes_{\textbf{k}[F]}\textbf{k}[H]\otimes_{\textbf{k}[G]}\langle [x_1,x_2]\rangle,$
\item [N79.] $\langle [x_2^-,[x_2^-,[x_2^-,x_1^-]]] \rangle \otimes_{\textbf{k}[F]}\textbf{k}[H]\otimes_{\textbf{k}[G]}\langle [x_1,x_2]\rangle,$
\item [N80.] $\langle [[x_1^-,x_2^-],[x_2^-,[x_2^-,x_1^-]]] \rangle \otimes_{\textbf{k}[F]}\textbf{k}[H]\otimes_{\textbf{k}[G]}\langle [x_1,x_2]\rangle.$
\end{description}

For the last 4 cases we are going to use Proposition \ref{pbwU}. From \textbf{(A46)},
$[[[x_1,x_2],[x_2,[x_2,x_1]]],x_2^-]\sim [x_2,x_1]^2+\alpha_1 g_2f_2x_2[x_2,x_1]x_1+\alpha_2 g_2f_2[x_2,[x_2,x_1]]x_1+\alpha_3 g_2f_2[x_2,x_1]^2$. Since $x_2 \notin \langle [[x_1,x_2],[x_2,[x_2,x_1]]]\rangle$, using the proposition, we exclude:
\begin{description}
\item [N81.] $\langle x_2^- \rangle \otimes_{\textbf{k}[F]}\textbf{k}[H]\otimes_{\textbf{k}[G]}\langle [[x_1,x_2],[x_2,[x_2,x_1]]]\rangle,$
\item [N82.] $\langle [[[x_1^-,x_2^-],x_2^-],x_2^-] \rangle \otimes_{\textbf{k}[F]}\textbf{k}[H]\otimes_{\textbf{k}[G]}\langle [[x_1,x_2],[x_2,[x_2,x_1]]]\rangle.$
\end{description}

Analogously, from \textbf{(A16)} we have $[x_2,[[x_1^-,x_2^-],[x_2^-,[x_2^-,x_1^-]]]]\sim [x_2^-,x_1^-]^2+\alpha_1 [x_2^-,x_1^-]^2g_2f_2+\alpha_2 [x_2^-,[x_2^-,x_1^-]]x_1^-g_2f_2+\alpha_3 x_2^-[x_2^-,x_1^-]x_1^-g_2f_2$. As we know that $x_2^- \notin \langle [[x_1^-,x_2^-],[x_2^-,[x_2^-,x_1^-]]]\rangle$, again from Proposition \ref{pbwU} we eliminate the final possibilities:
\begin{description}
\item [N83.] $\langle [[x_1^-,x_2^-],[x_2^-,[x_2^-,x_1^-]]] \rangle \otimes_{\textbf{k}[F]}\textbf{k}[H]\otimes_{\textbf{k}[G]}\langle x_2\rangle,$
\item [N84.] $\langle [[x_1^-,x_2^-],[x_2^-,[x_2^-,x_1^-]]] \rangle \otimes_{\textbf{k}[F]}\textbf{k}[H]\otimes_{\textbf{k}[G]}\langle [[[x_1,x_2],x_2],x_2]\rangle.$
\end{description}
Now we have showed that the first 60 cases are the only (homogeneous) right coideal subalgebras of $U_q(\mathfrak{g})$ (respectively, $u_q(\mathfrak{g})$) that contain $\textbf{k}[H]$. The theorem is proved. \hfill $\square$ \vspace{0.5cm}

\section*{Appendix}

\quad In this appendix we are going to calculate the skew commutators that are listed in the previous section. Below by $u,v,w$ we
denote homogeneous linear combinations of words in $x_1,x_2,\ldots x_n,$ while
$u^-, v^-, w^-$ are the elements that appeared from $u,v,w$ by replacing all $x_i$
with $x_i^-.$

We notice that, from Definition \ref{G2}, we have:
$$p(u,v)=p_{uv},\quad p(u,v^-)=p_{vu}, \quad p(u^-,v)=p_{uv}^{-1}, \quad p(u^-,v^-)=p_{vu}^{-1}.$$

We also have the following relations, that are a consequence of \eqref{gfx} and \eqref{gfx2}:
$$x_1g_1=p_{11}g_1x_1,\quad x_1g_2=p_{12}g_2x_1,\quad
x_1f_1=p_{11}f_1x_1,\quad x_1f_2=p_{21}f_2x_1,$$
$$x_2g_1=p_{21}g_1x_2,\quad x_2g_2=p_{22}g_2x_2,\quad
x_2f_1=p_{12}f_1x_2,\quad x_2f_2=p_{22}f_2x_2,$$
$$g_1x_1^-=p_{11}x_1^-g_1,\quad
g_2x_1^-=p_{12}x_1^-g_2,\quad f_1x_1^-=p_{11}x_1^-f_1,\quad
f_2x_1^-=p_{21}x_1^-f_2,$$
$$g_1x_2^-=p_{21}x_2^-g_1,\quad g_2x_2^-=p_{22}x_2^-g_2,\quad
f_1x_2^-=p_{12}x_2^-f_1,\quad f_2x_2^-=p_{22}x_2^-f_2.$$
\noindent where $p_{22}=q,$ $p_{11}
=q^3,$ $p_{12}p_{21}=q^{-3}.$

For simplification, we are also going to use the formulas, see \cite{Kh6} page 2612:

\begin{equation*}
[[u,v],w^-]=[u,[v,w^-]]+p_{wv}[[u,w^-],v],
\end{equation*}

\begin{equation*}
[u,[v^-,w^-]]=[[u,v^-],w^-]+p_{vu}[v^-,[u,w^-]],
\end{equation*}

\begin{equation*}
[u\cdot v,w]=p_{vw}[u,w]\cdot v+u\cdot [v,w],
\end{equation*}

\begin{equation*}
[u,v\cdot w]=[u,v]\cdot w+p_{uv}v\cdot [u,w].
\end{equation*}

Now we are ready to calculate the skew-commutators.

\begin{description}
\item [(A1)] $[x_1,x_1^-]=1-g_1f_1$,

\item [(A2)] $[x_1,x_2^-]=0$,

\item [(A3)]
$[x_1,[x_1^-,x_2^-]]=[[x_1,x_1^-],x_2^-]+p_{11}[x_1^-,[x_1,x_2^-]]=(1-g_1f_1)x_2^--x_2^-(1-g_1f_1)=-g_1f_1x_2^-+x_2^-g_1f_1=-p_{12}p_{21}x_2^-g_1f_1+x_2^-g_1f_1=(1-q^{-3})x_2^-g_1f_1$,

\item [(A4)]
$[x_1,[x_2^-,x_1^-]]=[x_1,x_2^-x_1^-]-p_{12}^{-1}[x_1,x_1^-x_2^-]=[x_1,x_2^-]x_1^-+p_{21}x_2^-[x_1,x_1^-]-p_{12}^{-1}[x_1,x_1^-]x_2^--p_{12}^{-1}p_{11}x_1^-[x_1,x_2^-]=p_{21}x_2^-(1-g_1f_1)-p_{12}^{-1}(1-g_1f_1)x_2^-=p_{21}x_2^--p_{21}x_2^-g_1f_1-p_{12}^{-1}x_2^-+p_{12}^{-1}p_{12}p_{21}x_2^-g_1f_1=(p_{21}-p_{12}^{-1})x_2^-=p_{21}(1-q^3)x_2^-$,

\item [(A5)]
$[x_1,[[x_1^-,x_2^-],x_2^-]]=[[x_1,[x_1^-,x_2^-]],x_2^-]+p_{11}p_{21}[[x_1^-,x_2^-],[x_1,x_2^-]]=(1-q^{-3})x_2^-g_1f_1x_2^--p_{22}^{-1}x_2^-(1-q^{-3})x_2^-g_1f_1=(1-q^{-3})(p_{12}p_{21}-p_{22}^{-1})(x_2^-)^2g_1f_1=-q^{-1}(1-q^{-2})(1-q^{-3})(x_2^-)^2g_1f_1$,

\item [(A6)]
$[x_1,[[[x_1^-,x_2^-],x_2^-],x_2^-]]=[[x_1,[[x_1^-,x_2^-],x_2^-]],x_2^-]+p_{11}p_{21}^2[[[x_1^-,x_2^-],x_2^-],[x_1,x_2^-]]=-q^{-1}(1-q^{-3})(1-q^{-2})(x_2^-)^2g_1f_1x_2^-+q^{-2}q^{-1}(1-q^{-3})(1-q^{-2})(x_2^-)^3g_1f_1=q^{-3}(1-q^{-1})(1-q^{-2})(1-q^{-3})(x_2^-)^3g_1f_1$,

\item [(A7)]
$[x_1,[[x_1^-,x_2^-],[[x_1^-,x_2^-],x_2^-]]]=[[x_1,[x_1^-,x_2^-]],[[x_1^-,x_2^-],x_2^-]]+p_{11}p_{21}[[x_1^-,x_2^-],[x_1,[[x_1^-,x_2^-],x_2^-]]]=(1-q^{-3})x_2^-g_1f_1[[x_1^-,x_2^-],x_2^-]-(1-q^{-3})p_{11}p_{21}^2p_{11}^{-1}p_{21}^{-2}p_{12}^{-1}p_{22}^{-2}[[x_1^-,x_2^-],x_2^-]x_2^-g_1f_1-p_{11}p_{21}q^{-1}(1-q^{-3})(1-q^{-2})[x_1^-,x_2^-](x_2^-)^2g_1f_1+p_{11}p_{21}p_{11}^{-1}p_{21}^{-1}p_{21}^{-1}p_{22}^{-2}q^{-1}(1-q^{-3})(1-q^{-2})(x_2^-)^2g_1f_1[x_1^-,x_2^-]=(1-q^{-3})q^6q^{-6}x_2^-[[x_1^-,x_2^-],x_2^-]g_1f_1-p_{21}q(1-q^{-3})([[[x_1^-,x_2^-],x_2^-],x_2^-]+p_{21}^{-1}q^{-2}x_2^-[[x_1^-,x_2^-],x_2^-])g_1f_1-p_{21}q^2(1-q^{-3})(1-q^{-2})([[[x_1^-,x_2^-],x_2^-],x_2^-]+p_{12}q(1+q)x_2^-[[x_1^-,x_2^-],x_2^-]+p_{12}^2q^4(x_2^-)^2[x_1^-,x_2^-])g_1f_1+p_{12}(1-q^{-3})(1-q^{-2})q^6q^{-3}(x_2^-)^2[x_1^-,x_2^-]g_1f_1=(1-q^{-3})(1-q^{-1}-(1-q^{-2})(1+q))x_2^-[[x_1^-,x_2^-],x_2^-]g_1f_1-p_{21}q(1-q^{-3})(1+q(1-q^{-2}))[[[x_1^-,x_2^-],x_2^-],x_2^-]g_1f_1+p_{12}(1-q^{-3})(1-q^{-2})(q^3-q^{3})(x_2^-)^2[x_1^-,x_2^-]g_1f_1=-q(1-q^{-3})^2x_2^-[[x_1^-,x_2^-],x_2^-]g_1f_1-p_{21}q(1-q^{-3})(1+q-q^{-1})[[[x_1^-,x_2^-],x_2^-],x_2^-]g_1f_1$,

\item [(A8)] $[x_2,x_1^-]=0$,

\item [(A9)] $[x_2,x_2^-]=1-g_2f_2$,

\item [(A10)] $[x_2,[x_1^-,x_2^-]]=p_{12}(1-q^3)x_1^-$,

\item [(A11)] $[x_2,[x_2^-,x_1^-]]=(1-q^{-3})x_1^-g_2f_2$,

\item [(A12)]
$[x_2,[[x_1^-,x_2^-],x_2^-]]=[x_2,[x_1^-,x_2^-]x_2^-]-p_{21}^{-1}p_{22}^{-1}[x_2,x_2^-[x_1^-,x_2^-]]=[x_2,[x_1^-,x_2^-]]x_2^-+p_{12}p_{22}[x_1^-,x_2^-][x_2,x_2^-]-p_{21}^{-1}p_{22}^{-1}[x_2,x_2^-][x_1^-,x_2^-]-p_{21}^{-1}p_{22}^{-1}p_{22}x_2^-[x_2,[x_1^-,x_2^-]]=p_{12}(1-q^3)x_1^-x_2^-+p_{12}p_{22}[x_1^-,x_2^-](1-g_2f_2)-p_{21}^{-1}p_{22}^{-1}(1-g_2f_2)[x_1^-,x_2^-]-p_{22}^{-1}p_{21}^{-1}p_{22}x_2^-p_{12}(1-q^3)x_1^-=p_{12}(1-q^3)[x_1^-,x_2^-]+p_{12}q[x_1^-,x_2^-]-p_{12}q^2[x_1^-,x_2^-]-p_{12}q[x_1^-,x_2^-]g_2f_2+p_{12}q^2g_2f_2[x_1^-,x_2^-]=(p_{12}q^4q^{-3}-p_{12}q)g_2f_2[x_1^-,x_2^-]+p_{12}(1+q-q^2-q^3)[x_1^-,x_2^-]=p_{12}(1+q)(1-q^2)[x_1^-,x_2^-]$,

\item [(A13)]
$[x_2,[x_2^-,[x_2^-,x_1^-]]]=[x_2,x_2^-[x_2^-,x_1^-]]-p_{22}^{-1}p_{12}^{-1}[x_2,[x_2^-,x_1^-]x_2^-]=[x_2,x_2^-][x_2^-,x_1^-]+p_{22}x_2^-[x_2,[x_2^-,x_1^-]]-p_{22}^{-1}p_{12}^{-1}[x_2,[x_2^-,x_1^-]]x_2^--p_{22}^{-1}p_{12}^{-1}p_{22}p_{12}[x_2^-,x_1^-][x_2,x_2^-]=(1-g_2f_2)[x_2^-,x_1^-]+p_{22}x_2^-(1-q^{-3})x_1^-g_2f_2-p_{22}^{-1}p_{12}^{-1}(1-q^{-3})x_1^-g_2f_2x_2^--[x_2^-,x_1^-](1-g_2f_2)=[x_2^-,x_1^-]-g_2f_2[x_2^-,x_1^-]+p_{22}(1-q^{-3})x_2^-x_1^-g_2f_2-p_{22}p_{12}^{-1}(1-q^{-3})x_1^-x_2^-g_2f_2-[x_2^-,x_1^-]+[x_2^-,x_1^-]g_2f_2=(1-q^{-1})[x_2^-,x_1^-]g_2f_2+q(1-q^{-3})[x_2^-,x_1^-]g_2f_2=(1+q)(1-q^{-2})[x_2^-,x_1^-]g_2f_2$,

\item [(A14)]
$[x_2,[[[x_1^-,x_2^-],x_2^-],x_2^-]]=[x_2,[[x_1^-,x_2^-],x_2^-]x_2^-]-p_{22}^{-2}p_{21}^{-1}[x_2,x_2^-[[x_1^-,x_2^-],x_2^-]]=[x_2,[[x_1^-,x_2^-],x_2^-]]x_2^-+p_{12}p_{22}^2[[x_1^-,x_2^-],x_2^-][x_2,x_2^-]-p_{22}^{-2}p_{21}^{-1}[x_2,x_2^-][[x_1^-,x_2^-],x_2^-]-p_{22}^{-2}p_{21}^{-1}p_{22}x_2^-[x_2,[[x_1^-,x_2^-],x_2^-]]=p_{12}(1+q-q^2-q^3)[x_1^-,x_2^-]x_2^-+p_{12}p_{22}^2[[x_1^-,x_2^-],x_2^-]$ $(1-g_2f_2)-p_{22}^{-2}p_{21}^{-1}(1-g_2f_2)[[x_1^-,x_2^-],x_2^-]-p_{21}^{-1}p_{22}^{-1}p_{12}(1+q-q^2-q^3)x_2^-[x_1^-,x_2^-]=p_{12}(1+q-q^2-q^3)[[x_1^-,x_2^-],x_2^-]+(p_{12}q^2-p_{12}q)[[x_1^-,x_2^-],x_2^-]-(p_{12}q^2-p_{12}q^2)[[x_1^-,x_2^-],x_2^-]g_2f_2=p_{12}(1-q^3)[[x_1^-,x_2^-],x_2^-]$,

\item [(A15)]
$[x_2,[x_2^-,[x_2^-,[x_2^-,x_1^-]]]]=[x_2,x_2^-[x_2^-,[x_2^-,x_1^-]]]-p_{22}^{-2}p_{12}^{-1}[x_2,[x_2^-,[x_2^-,x_1^-]]x_2^-]=[x_2,x_2^-][x_2^-,[x_2^-,x_1^-]]+p_{22}x_2^-[x_2,[x_2^-,[x_2^-,x_1^-]]]-p_{22}^{-2}p_{12}^{-1}[x_2,[x_2^-,[x_2^-,x_1^-]]]x_2^--p_{22}^{-2}p_{12}^{-1}p_{22}^2p_{12}[x_2^-,[x_2^-,x_1^-]][x_2,x_2^-]=(1-g_2f_2)[x_2^-,[x_2^-,x_1^-]]+p_{22}(1+q-q^{-1}-q^{-2})x_2^-[x_2^-,x_1^-]g_2f_2-p_{22}^{-2}p_{12}^{-1}(1+q-q^{-1}-q^{-2})[x_2^-,x_1^-]g_2f_2x_2^--[x_2^-,[x_2^-,x_1^-]](1-g_2f_2)=p_{22}(1+q-q^{-1}-q^{-2})[x_2^-,[x_2^-,x_1^-]]g_2f_2=q^2(1-q^{-3})[x_2^-,[x_2^-,x_1^-]]g_2f_2$,

\item [(A16)]
$[x_2,[[x_1^-,x_2^-],[x_2^-,[x_2^-,x_1^-]]]]=[[x_2,[x_1^-,x_2^-]],[x_2^-,[x_2^-,x_1^-]]]+p_{12}p_{22}[[x_1^-,x_2^-],$ $[x_2,[x_2^-,[x_2^-,x_1^-]]]]=p_{12}(1-q^3)x_1^-[x_2^-,[x_2^-,x_1^-]]-p_{22}^2p_{12}p_{21}^{-2}p_{22}^{-2}p_{11}^{-1}p_{12}^{-1}p_{12}(1-q^3)[x_2^-,[x_2^-,x_1^-]]x_1^-+p_{12}p_{22}(1+q-q^{-1}-q^{-2})[x_1^-,x_2^-][x_2^-,x_1^-]g_2f_2-p_{12}^{-1}p_{22}^{-1}p_{21}^{-2}p_{22}^{-2}p_{11}^{-1}p_{12}^{-1}p_{12}p_{22}(1+q-q^{-1}-q^{-2})[x_2^-,x_1^-]g_2f_2[x_1^-,x_2^-]=p_{12}(1-q^3)(p_{12}q^2(1-q^{-2})[x_2^-,x_1^-]^2+p_{12}^2q^3[x_2^-,[x_2^-,x_1^-]]x_1^-)-p_{12}^3q^3(1-q^3)[x_2^-,[x_2^-,x_1^-]]x_1^-+p_{12}q(1+q-q^{-1}-q^{-2})(p_{12}(1-q^3)x_2^-x_1^--p_{12}[x_2^-,x_1^-])[x_2^-,x_1^-]g_2f_2-p_{12}(1+q-q^{-1}-q^{-2})[x_2^-,x_1^-]q^{-1}(p_{12}(1-q^3)x_2^-x_1^--p_{12}[x_2^-,x_1^-])g_2f_2=p_{12}^2q^2(1-q^3)(1-q^{-2})[x_2^-,x_1^-]^2-p_{12}^2q(1+q)(1-q^{-2})^2[x_2^-,x_1^-]^2g_2f_2-p_{12}^3(1+q)(1-q^3)(1-q^{-2})[x_2^-,[x_2^-,x_1^-]]x_1^-g_2f_2-p_{12}^3(1+q)(1-q^{3})(1-q^{4})(1-q^{-2})x_2^-[x_2^-,x_1^-]x_1^-g_2f_2$,

\item [(A17)]
$[[x_1,x_2],x_1^-]=[x_1x_2,x_1^-]-p_{12}[x_2x_1,x_1^-]=p_{12}[x_1,x_1^-]x_2+x_1[x_2,x_1^-]-p_{12}p_{11}[x_2,x_1^-]x_1-p_{12}x_2[x_1,x_1^-]=p_{12}(1-g_1f_1)x_2-p_{12}x_2(1-g_1f_1)=p_{12}^2p_{21}g_1f_1x_2-p_{12}g_1f_1x_2=-p_{12}(1-q^{-3})g_1f_1x_2$,

\item [(A18)]
$[[x_1,x_2],x_2^-]=[x_1x_2,x_2^-]-p_{12}[x_2x_1,x_2^-]=p_{22}[x_1,x_2^-]x_2+x_1[x_2,x_2^-]-p_{12}p_{21}[x_2,x_2^-]x_1-p_{12}x_2[x_1,x_2^-]=x_1(1-g_2f_2)-p_{12}p_{21}(1-g_2f_2)x_1=(1-q^{-3})x_1$,

\item [(A19)]
$[[x_1,x_2],[x_2^-,x_1^-]]=[x_1x_2,[x_2^-,x_1^-]]-p_{12}[x_2x_1,[x_2^-,x_1^-]]=p_{12}p_{22}[x_1,[x_2^-,x_1^-]]x_2+x_1[x_2,[x_2^-,x_1^-]]-p_{12}p_{11}p_{21}[x_2,[x_2^-,x_1^-]]x_1-p_{12}x_2[x_1,[x_2^-,x_1^-]]=p_{12}p_{21}q(1-q^{3})x_2^-x_2+(1-q^{-3})x_1x_1^-g_2f_2-(1-q^{-3})x_1^-g_2f_2x_1-p_{12}p_{21}q(1-q^{3})x_2x_2^-=(q^3-1)x_1^-x_1g_2f_2+(1-q^{-3})g_2f_2-(1-q^{-3})g_1f_1g_2f_2-(1-q^{-3})x_1^-g_2f_2x_1-q^{-3}(1-q^3)+q^{-3}(1-q^3)g_2f_2=(1-q^{-3})(1-g_1f_1g_2f_2)$,

\item [(A20)]
$[[x_1,x_2],[x^-_2,[x^-_2,x^-_1]]]=[[[x_1,x_2],x^-_2],[x^-_2,x^-_1]]+p_{21}p_{22}[x^-_2,[[x_1,x_2],[x^-_2,x^-_1]]]=(1-q^{-3})x_1[x^-_2,x^-_1]-p_{21}p_{22}p_{11}p_{12}p_{22}^{-1}p_{12}^{-1}(1-q^{-3})[x^-_2,x^-_1]x_1+p_{21}p_{22}(1-q^{-3})x^-_2(1-g_1f_1g_2f_2)-p_{21}p_{22}(1-q^{-3})p_{21}^{-1}p_{22}^{-1}p_{22}^{-1}p_{12}^{-1}(1-g_1f_1g_2f_2)x^-_2=(1-q^{-3})[x_1,[x^-_2,x^-_1]]+p_{21}p_{22}(1-q^{-3})(1-q)x^-_2=p_{21}(1-q^{-3})(1-q^{-3})x^-_2+p_{21}p_{22}q(1-q^{-3})(1-q)x^-_2=p_{21}(1+q)(1-q^{2})(1-q^{-3})x^-_2$,

\item [(A21)] $[[x_2,x_1],x_1^-]=(1-q^{-3})x_2$,

\item [(A22)] $[[x_2,x_1],x_2^-]=-p_{21}(1-q^{-3})g_2f_2x_1$,

\item [(A23)] $[[x_2,x_1],[x_1^-,x_2^-]]=(1-q^{-3})(1-g_2f_2g_1f_1)$,

\item [(A24)] $[[x_2,x_1],[x^-_1,x^-_2],x^-_2]=[[[x_2,x_1],[x^-_1,x^-_2]],x^-_2]+p_{12}p_{11}p_{22}p_{21}[[x^-_1,x^-_2],[[x_2,x_1],x^-_2]]=(1-q^{-3})(1-g_2f_2g_1f_1)x^-_2-p_{22}p_{21}p_{21}^{-1}p_{22}^{-1}(1-q^{-3})x^-_2(1-g_2f_2g_1f_1)-qp_{21}(1-q^{-3})[x^-_1,x^-_2]g_2f_2x_1+qp_{21}(1-q^{-3})p_{21}^{-1}p_{22}^{-1}p_{12}^{-1}p_{11}^{-1}p_{22}^{-1}p_{21}^{-1}g_2f_2x_1[x^-_1,x^-_2]=(1-q^{-3})(1-q^{-1})x^-_2g_2f_2g_1f_1-p_{21}q(1-q^{-3})[x^-_1,x^-_2]g_2f_2x_1+q^{-1}(1-q^{-3})g_2f_2((1-q^{-3})x^-_2g_1f_1+p_{11}p_{21}[x^-_1,x^-_2]x_1)=(1-q^{-3})(1-q^{-1})x^-_2g_2f_2g_1f_1-p_{21}q(1-q^{-3})[x^-_1,x^-_2]g_2f_2x_1+q(1-q^{-3})^2x^-_2g_2f_2g_1f_1+p_{21}q^2(1-q^{-3})p_{12}p_{21}q^2$ $[x^-_1,x^-_2]g_2f_2x_1=(1-q)(1-q^{-2})(1-q^{-3})x^-_2g_1f_1g_2f_2$,

\item [(A25)] $[[x_2,x_1],[[[x^-_1,x^-_2],x^-_2],x^-_2]]=
[[[x_2,x_1],[[x^-_1,x^-_2],x^-_2]],x^-_2]+p_{12}p_{11}p_{22}^2p_{21}^2$ $[[[x^-_1,x^-_2],x^-_2],[[x_2,x_1],x^-_2]]=
(1-q^{-3})(1-q^{-2})(1+q)x^-_2g_1f_1g_2f_2x^-_2-p_{22}p_{21}p_{21}^{-1}p_{22}^{-2}(1-q^{-3})(1-q^{-2})(1+q)(x^-_2)^2g_1f_1g_2f_2-p_{21}^2q^2
[[x^-_1,x^-_2],x^-_2](1-q^{-3})g_2f_2x_1+p_{21}^2q^2(1-q^{-3})p_{12}^{-1}p_{11}^{-1}p_{22}^{-2}p_{21}^{-2}p_{21}^{-1}p_{22}^{-2}g_2f_2x_1[[x^-_1,x^-_2],x^-_2]=
-p_{21}^{2}q^2(1-q^{-3})[[x^-_1,x^-_2],x^-_2]g_2f_2x_1+q^{-2}(1-q^{-3})g_2f_2(p_{21}^2p_{11}[[x^-_1,x^-_2],x^-_2]x_1-q^{-1}(1-q^{-3})(1-q^{-2})(x^-_2)^2g_1f_1)=
-q(1-q^{-2})(1-q^{-3})^2(x^-_2)^2g_1f_1g_2f_2$,

\item [(A26)] $[[x_2,x_1],[[x^-_1,x^-_2],[[x^-_1,x^-_2],x^-_2]]]=
[[[x_2,x_1],[x^-_1,x^-_2]],[[x^-_1,x^-_2],x^-_2]]+p_{12}p_{11}\times p_{22}p_{21}[[x^-_2,x^-_1],[[x_1,x_2],[[x^-_1,x^-_2],x^-_2]]]=
(1-q^{-3})(1-g_2f_2g_1f_1)[[x^-_1,x^-_2],x^-_2]-p_{12}p_{11}p_{22}^2p_{21}^2p_{11}^{-1}p_{12}^{-1}p_{21}^{-2}p_{22}^{-2}(1-q^{-3})[[x^-_1,x^-_2],x^-_2]
(1-g_2f_2g_1f_1)+p_{22}(1-q^{-3})(1-q^{-2})(1+q)[x^-_1,x^-_2]x^-_2g_2f_2g_1f_1
-p_{22}(1-q^{-3})(1-q^{-2})(1+q)p_{12}^{-2}p_{11}^{-2}p_{22}^{-2}p_{21}^{-2}\times x^-_2g_2f_2g_1f_1[x^-_1,x^-_2]=
(1-q)(1-q^{-3})[[x^-_1,x^-_2],x^-_2]$ $g_1f_1g_2f_2+p_{22}(1-q^{-3})(1-q^{-2})(1+q)[x^-_1,x^-_2]x^-_2g_2f_2g_1f_1
-p_{12}p_{22}^3(1-q^{-3})(1-q^{-2})(1+q)x^-_2[x^-_1,x^-_2]g_2f_2g_1f_1=
(1-q)(1-q^{-3})[[x^-_1,x^-_2],x^-_2]g_1f_1g_2f_2+p_{22}(1-q^{-3})(1-q^{-2})(1+q)[[x^-_1,x^-_2],x^-_2]g_2f_1g_1f_1=
(1-q^{-3})(1-q+q(1-q^{-2})(1+q))[[x^-_1,x^-_2],x^-_2]g_2f_2g_1f_1=
q^2(1-q^{-3})[[x^-_1,x^-_2],x^-_2]g_1f_1g_2f_2$,

\item [(A27)]
$[[[x_1,x_2],x_2],x_1^-]=[[x_1,x_2]x_2,x_1^-]-p_{12}p_{22}[x_2[x_1,x_2],x_1^-]=p_{12}[[x_1,x_2],x_1^-]x_2+[x_1,x_2][x_2,x_1^-]-p_{12}p_{22}p_{11}p_{12}[x_2,x_1^-][x_1,x_2]-p_{12}p_{22}x_2[[x_1,x_2],x_1^-]=-p_{12}^2(1-q^{-3})g_1f_1x_2^2+p_{12}^3p_{21}p_{22}(1-q^{-3})g_1f_1x_2^2=-p_{12}^2(1-q^{-2})(1-q^{-3})g_1f_1x_2^2$,

\item [(A28)]
$[[[x_1,x_2],x_2],x_2^-]=[[x_1,x_2]x_2,x_2^-]-p_{12}p_{22}[x_2[x_1,x_2],x_2^-]=p_{22}[[x_1,x_2],x_2^-]x_2+[x_1,x_2][x_2,x_2^-]-p_{12}p_{22}p_{21}p_{22}[x_2,x_2^-][x_1,x_2]-p_{12}p_{22}x_2[[x_1,x_2],x_2^-]=p_{22}(1-q^{-3})x_1x_2+[x_1,x_2](1-g_2f_2)-q^{-1}(1-g_2f_2)-p_{12}p_{22}(1-q^{-3})x_2x_1=(1+q)(1-q^{-2})[x_1,x_2]$,

\item [(A29)]
$[[[x_1,x_2],x_2],[x^-_2,x^-_1]]=
[[x_1,x_2],[x_2,[x^-_2,x^-_1]]]+p_{22}p_{12}[[[x_1,x_2],[x^-_2,x^-_1]],x_2]=
(1-q^{-3})[x_1,x_2]x^-_1g_2f_2-p_{12}p_{22}p_{21}p_{22}p_{11}p_{12}(1-q^{-3})x^-_1g_2f_2[x_1,x_2]
+p_{12}p_{22}(1-q^{-3})(1-g_1f_1g_2f_2)x_2-p_{12}p_{22}p_{12}p_{22}p_{22}^{-1}p_{12}^{-1}(1-q^{-3})x_2(1-g_1f_1g_2f_2)=
(1-q^{-3})(p_{12}p_{11}x^-_1[x_1,x_2]-p_{12}(1-q^{-3})g_1f_1x_2)g_2f_2-p_{12}q^2(1-q^{-3})x^-_1g_2f_2[x_1,x_2]
-(1-q^{-1})p_{12}q(1-q^{-3})g_1f_1g_2f_2x_2=
-p_{12}(1-q^{-3})^2q^2g_1f_1g_2f_2x_2-p_{12}q(1-q^{-1})(1-q^{-3})g_1f_1g_2f_2x_2=
-p_{12}q(1+q)(1-q^{-2})(1-q^{-3})g_1f_1g_2f_2x_2$,

\item [(A30)]
$[[[x_1,x_2],x_2],[x^-_2,[x^-_2,x^-_1]]]=
[[[x_1,x_2],x_2],[x^-_2,x^-_1]]+p_{21}q^2[x^-_2,[[x_1,x_2],x_2]]=
(1+q-q^{-1}-q^{-2})[x_1,x_2][x^-_2,x^-_1]-p_{22}^{-1}p_{12}^{-1}p_{21}p_{22}^2p_{11}p_{12}^2(1+q-q^{-1}-q^{-2})[x^-_2,x^-_1][x_1,x_2]
-p_{21}q^2p_{12}q(1-q^{-3})(1-q^{-2})(1+q)x^-_2g_1f_1g_2f_2x_2+(1-q^{-3})(1-q^{-2})(1+q)p_{22}^{-1}p_{12}^{-1}p_{21}p_{22}^2p_{11}p_{12}^2g_1f_1g_2f_2x_2x^-_2=
(1+q-q^{-1}-q^{-2})((1-q^{-3})(1-g_1f_1g_2f_2)+q[x^-_2,x^-_1][x_1,x_2])-q(1+q-q^{-1}-q^{-2})[x^-_2,x^-_1][x_1,x_2]
-(1-q^{-3})(1-q^{-2})(1+q)x^-_2g_1f_1g_2f_2x_2+(1-q^{-3})(1-q^{-2})(1+q)g_1f_1g_2f_2(1-g_2f_2+qx^-_2x_2)=
(1+q-q^{-1}-q^{-2})(1-q^{-3})(1-g_1f_1g_2f_2)-(1-q^{-3})(1-q^{-2})(1+q)x^-_2g_1f_1g_2f_2x_2
+(1-q^{-3})(1-q^{-2})(1+q)g_1f_1g_2f_2(1-g_2f_2)+(1-q^{-3})(1-q^{-2})(1+q)x^-_2g_1f_1g_2f_2x_2=
(1+q)(1-q^{-2})(1-q^{-3})(1-g_1f_1g_2^2f_2^2)$,

\item [(A31)]
$[[[x_1,x_2],x_2],[x_{2}^{-},[x_{2}^{-},[x_{2}^{-},x_{1}^{-}]]]]=[[[x_1,x_2],x_2],x_{2}^{-}],[x_{2}^{-},[x_{2}^{-},x_{1}^{-}]]
 + p_{21}q^{2}$ $[x_{2}^{-},[[[x_1 ,x_2],x_2],[x_{2}^{-},[x_{2}^{-},x_{1}^{-}]]]]
 = (1+q-q^{-1}-q^{-2})[x_1,x_2][x_{2}^{-},[x_{2}^{-},x_{1}^{-}]]
 - p_{21}^{2}p_{22}^{4}p_{11}p_{12}^{2}p_{22}^{-2}p_{12}^{-1}(1+q-q^{-1}-q^{-2})[x_{2}^{-},[x_{2}^{-},x_{1}^{-}]][x_1,x_2]
 + p_{21}q^{2}(1-q^{-3})(1-q^{-2})(1+q)x_{2}^{-}(1-g_{1}f_{1}g_{2}^2f_{2}^2)
 - p_{21}q^2p_{21}^{-1}p_{22}^{-2}p_{22}^{-2}p_{12}^{-1}(1-q^{-3})(1-q^{-2})(1+q)(1-g_{1}f_{1}g_{2}^2f_{2}^2)x_{2}^{-}
 = (1+q-q^{-1}-q^{-2})(p_{21}(1-q^{-3})(1-q^{-2})(1+q)x_{2}^{-}
 + p_{21}^{2}p_{22}^{2}p_{11}p_{12}[x_{2}^{-},[x_{2}^{-},x_{1}^{-}]][x_1,x_2])
 - p_{21}q^2(1+q-q^{-1}-q^{-2})[x_{2}^{-},[x_{2}^{-},x_{1}^{-}]][x_1,x_2]
 + p_{21}q^2(1-q^{-1})(1-q^{-3})(1-q^{-2})(1+q)x_{2}^{-}
 = p_{21}(1-q^{-3})(1+q)((1-q^{-2})(1+q-q^{-1}-q^{-2})
 + q^2(1-q^{-1})(1-q^{-2}))x_{2}^{-}= p_{21}q(1+q)(1-q^{-2})(1-q^{-3})^{2}x_{2}^{-}$,

\item [(A32)]
$[[x_2,[x_2,x_1]],x_1^-]=[x_2[x_2,x_1],x_1^-]-p_{22}p_{21}[[x_2,x_1]x_2,x_1^-]=p_{12}p_{11}[x_2,x_1^-][x_2,x_1]+x_2[[x_2,x_1],x_1^-]-p_{22}p_{12}p_{21}[[x_2,x_1],x_1^-]x_2-p_{22}p_{21}[x_2,x_1][x_2,x_1^-]=(1-q^{-3})x_2^2-q^{-2}(1-q^{-3})x_2^2=(1-q^{-2})(1-q^{-3})x_2^2$,

\item [(A33)]
$[[x_2,[x_2,x_1]],x_2^-]=[x_2[x_2,x_1],x_2^-]-p_{22}p_{21}[[x_2,x_1]x_2,x_2^-]=p_{22}p_{21}[x_2,x_2^-][x_2,x_1]+x_2[[x_2,x_1],x_2^-]-p_{22}p_{22}p_{21}[[x_2,x_1],x_2^-]x_2-p_{22}p_{21}[x_2,x_1][x_2,x_2^-]=p_{22}p_{21}(1-g_2f_2)[x_2,x_1]-p_{21}(1-q^{-3})x_2g_2f_2x_1+p_{22}^2p_{21}^2(1-q^{-3})g_2f_2x_1x_2-p_{22}p_{21}[x_2,x_1](1-g_2f_2)=-p_{21}q^2(1-q^{-3})g_2f_2[x_2,x_1]-p_{21}qg_2f_2[x_2,x_1]+p_{21}g_2f_2[x_2,x_1]=-p_{21}q(1+q)(1-q^{-2})g_2f_2[x_2,x_1]$,

\item [(A34)]
$[[x_2,[x_2,x_1]],[x_{1}^{-},x_{2}^{-}]]= [x_2,[[x_2,x_1],[x_{1}^{-},x_{2}^{-}]]]
+ p_{12}p_{11}p_{22}p_{21}[[x_2,[x_{1}^{-},x_{2}^{-}]],[x_2,x_1]]
= (1-q^{-3})x_2(1-g_2f_2g_1f_1)-p_{22}p_{21}p_{12}p_{22}(1-q^{-3})(1-g_1f_1g_2f_2)x_2
+ p_{12}(1-q^{3})x_{1}^{-}[x_2,x_1]-p_{22}p_{21}p_{12}^{-1}p_{11}^{-1}p_{22}^{-1}p_{21}^{-1}p_{12}(1-q^{3})[x_2,x_1]x_{2}^{-}
= (1-q^{-3})(1-q^{-1})x_2 + p_{12}(1-q^{3})x_{1}^{-}[x_2,x_1]
- q^{-3}(1-q^{3})((1-q^{-3})x_2 + p_{12}q^3x_{1}^{-}[x_2,x_1])
= (1+q)(1-q^{-2})(1-q^{-3})x_2,
$

\item [(A35)]
$
[[x_2,[x_2,x_1]],[[x_{1}^{-},x_{2}^{-}],x_{2}^{-}]]= [x_2,[[x_2,x_1],[[x_{1}^{-},x_{2}^{-}],x_{2}^{-}]]]
+ p_{12}p_{22}^{2}p_{11}p_{21}^{2}$ $[[x_2,[[x_{1}^{-},x_{2}^{-}],x_{2}^{-}]],[x_2,x_1]]
= (1-q^{-3})(1-q^{-2})(1+q)x_2x_{2}^{-}g_2f_2g_1f_1
- p_{22}p_{21}p_{12}p_{22}^{2}(1-q^{-3})(1-q^{-2})(1+q)x_{2}^{-}g_2f_2g_1f_1x_2
+ p_{21}q^2p_{12}(1+q-q^2-q^3)[x_{1}^{-},x_{2}^{-}][x_2,x_1]
- q^{-1}p_{22}p_{21}p_{12}^{-1}p_{11}^{-1}p_{22}^{-2}p_{21}^{-2}(1+q-q^2-q^3)[x_2,x_1][x_{1}^{-},x_{2}^{-}]
= (1-q^{-3})(1-q^{-2})(1+q)(1-g_2f_2+qx_{2}^{-}x_2)g_2f_2g_1f_1
- (1-q^{-3})(1-q^{-2})(1+q)x_{2}^{-}g_2f_2g_1f_1x_2
+ q^{-1}(1+q-q^2-q^3)[x_{1}^{-},x_{2}^{-}][x_2,x_1]
- q^{-2}(1+q-q^2-q^3)((1-q^{-3})(1-g_2f_2g_1f_1)
+ p_{12}p_{11}p_{22}p_{21}[x_{1}^{-},x_{2}^{-}][x_2,x_1])
= (1+q)(1-q^{-2})(1-q^{-3})(1-g_1f_1g_{2}^{2}f_{2}^{2})
$,

\item [(A36)]
$
[[x_2,[x_2,x_1]],[[[x_{1}^{-},x_{2}^{-}],x_{2}^{-}],x_{2}^{-}]]= [[[x_2,[x_2,x_1]],[[x_{1}^{-},x_{2}^{-}],x_{2}^{-}]],x_{2}^{-}]
+ p_{12}^2$ $p_{11}p_{22}^4p_{21}^2[[[x_{1}^{-},x_{2}^{-}],x_{2}^{-}],[[x_2,[x_2,x_1]],x_{2}^{-}]]
= (1-q^{-3})(1-q^{-2})(1+q)(1-g_{2}^2f_{2}^2g_1f_1)x_{2}^- -p_{22}^2p_{21}p_{21}^{-1}p_{22}^{-2}(1-q^{-3})(1-q^{-2})(1+q)x_{2}^{-}(1-g_{2}^2f_{2}^2g_1f_1) + qp_{21}(1-q+q^{-1}-q^{2})[[x_{1}^{-},x_{2}^{-}],x_{2}^{-}]g_{2}f_{2}[x_2,x_1] - p_{12}^{-2}p_{11}^{-1}p_{22}^{-4}p_{21}^{-2}p_{21}^{-1}p_{22}^{-2}q
p_{21}(1-q+q^{-1}-q^{2})g_{2}f_{2}[x_2,x_1][[x_{1}^{-},x_{2}^{-}],x_{2}^{-}]=(1-q^{-3})(1-q^{-2})(1+q)(1-q)x_{2}^{-}g_{2}^2f_{2}^2g_1f_1 + p_{21}q(1-q+q^{-1}-q^{2})[[x_{1}^{-},x_{2}^{-}],x_{2}^{-}]g_{2}f_{2}[x_2,x_1] - q^{-2}(1-q+q^{-1}-q^{2})g_{2}f_{2}((1-q^{-3})(1-q^{-2})(1+q)x_{2}^{-}g_{2}^2f_{2}^2g_1f_1 + p_{12}p_{11}p_{22}^2p_{21}^2[[x_{1}^{-},x_{2}^{-}],x_{2}^{-}][x_2,x_1])= (1-q^{-3})(1-q^{-2})(1+q)(1-q-1+q+q^{-1}-q^{-2})x_{2}^{-}g_{2}^2f_{2}^2g_1f_1= q^2(1+q)(1-q^{-2})(1-q^{-3})^2x_{2}^{-}g_1f_1g_{2}^2f_{2}^2
$,

\item[(A37)]
$[[x_2,[x_2,x_1]],[[x_1^-,x_2^-],[[x_1^-,x_2^-],x_2^-]]]=-p_{12}q^3(1+q)(1-q^{-2})(1-q^{-3})^2[x_1^-,x_2^-]$,

\item [(A38)]
$[[[[x_1,x_2],x_2],x_2],x_1^-]=[[[x_1,x_2],x_2]x_2,x_1^-]-p_{12}p_{22}^2[x_2[[x_1,x_2],x_2],x_1^-]=p_{12}[[[x_1,x_2],x_2],x_1^-]x_2+[[x_1,x_2],x_2][x_2,x_1^-]-p_{12}p_{22}^2p_{11}p_{12}^2[x_2,x_1^-][[x_1,x_2],x_2]-p_{12}p_{22}^2x_2[[[x_1,x_2],x_2],x_1^-]=-p_{12}^3(1-q^{-3})(1-q^{-2})g_1f_1x_2^3+p_{12}^3q^2(1-q^{-3})(1-q^{-2})p_{12}p_{21}g_1f_1x_2^3=-p_{12}^3(1-q^{-1})(1-q^{-2})(1-q^{-3})g_1f_1x_2^3$,

\item [(A39)]
$[[[[x_1,x_2],x_2],x_2],x_2^-]=[[[x_1,x_2],x_2]x_2,x_2^-]-p_{12}p_{22}^2[x_2[[x_1,x_2],x_2],x_2^-]=p_{22}[[[x_1,x_2],x_2],x_2^-]x_2+[[x_1,x_2],x_2][x_2,x_2^-]-p_{12}p_{22}^4p_{21}[x_2,x_2^-][[x_1,x_2],x_2]-p_{12}p_{22}^2x_2[[[x_1,x_2],x_2],x_2^-]=p_{22}(1+q-q^{-1}-q^{-2})[x_1,x_2]x_2+[[x_1,x_2],x_2](1-g_2f_2)-q(1-g_2f_2)[[x_1,x_2],x_2]-p_{12}q^2(1+q-q^{-1}-q^{-2})x_2[x_1,x_2]=q(1+q-q^{-1}-q^{-2})[[x_1,x_2],x_2]+(1-q)[[x_1,x_2],x_2]-[[x_1,x_2],x_2]g_2f_2+qg_2f_2[[x_1,x_2],x_2]=q^2(1-q^{-3})[[x_1,x_2],x_2]$,

\item [(A40)]
$
[[[[x_1,x_2],x_2],x_2],[x_{2}^{-},x_{1}^{-}]]=[[[x_1,x_2],x_2],[x_2,[x_{2}^{-},x_{1}^{-}]]]+ p_{22}p_{12}[[[[x_1,x_2],x_2],$ $[x_{2}^{-},x_{1}^{-}]],x_2]=(1-q^{-3})[[x_1,x_2],x_2]x_{1}^{-}g_2f_2 - p_{12}p_{22}^2p_{21}p_{22}^2p_{11}p_{12}^2(1-q^{-3})x_{1}^{-}g_2f_2[[x_1,x_2],x_2]-p_{12}^2q^2(1-q^{-3})(1-q^{-2})(1+q)g_1f_1g_2f_2x_{2}^{2}
+p_{12}p_{22}^2p_{22}^{-1}p_{12}^{-1}p_{12}^{2}q^2(1-q^{-3})(1-q^{-2})(1+q)x_2g_1f_1g_2f_2x_{2}=(1-q^{-3})(-p_{12}^{2}(1-q^{-3})
(1-q^{-2})g_1f_1x_{2}^{2}+p_{11}p_{12}^{2}x_{1}^{-}[[x_1,x_2],x_2])g_2f_2 - p_{12}^{2}q^4(1-q^{-3})x_{1}^{-}g_2f_2[[x_1,x_2],x_2]= -p_{12}^2q^4(1-q^{-2})(1-q^{-3})^2g_1f_1g_2f_2x^2_2
$,

\item [(A41)]
$
[[[[x_1,x_2],x_2],x_2],[x_{2}^{-},[x_{2}^{-},x_{1}^{-}]]]=[[[x_1,x_2],x_2],[x_{2},[x_{2}^{-},[x_{2}^{-},x_{1}^{-}]]]]+
p_{22}^2p_{12}$ $[[[[x_1,x_2],x_2],[x_{2}^{-},[x_{2}^{-},x_{1}^{-}]]],x_2]=(1+q-q^{-1}-q^{-2})[[x_1,x_2],x_2][x_{2}^{-},x_{1}^{-}]g_2f_2-
p_{12}p_{22}^2p_{21}^2p_{22}^4p_{11}p_{12}^2(1+q-q^{-1}-q^{-2})[x_{2}^{-},x_{1}^{-}]g_2f_2[[x_1,x_2],x_2] +
p_{12}q^2(1-q^{-3})(1-q^{-2})(1+q)(1-g_1f_1g_2^2f_2^2)x_2-p_{12}p_{22}^2p_{22}^{-2}p_{12}^{-1}p_{12}q^2(1-q^{-3})(1-q^{-2})(1+q)x_2(1-g_1f_1g_2^2f_2^2)=
(1+q-q^{-1}-q^{-2})(-p_{12}q(1-q^{-3})(1-q^{-2})(1+q)g_1f_1g_2f_2x_2+p_{21}p_{22}^2p_{11}p_{12}^2[x_{2}^{-},x_{1}^{-}][[x_1,x_2],x_2])g_2f_2-
p_{12}q^3(1-q-q^{-1}-q^{-2})[x_{2}^{-},x_{1}^{-}]g_2f_2[[x_1,x_2],x_2]- p_{12}q^2(1-q^{-3})(1-q^{-2})(1+q)g_1f_1g_2^2f_2^2x_2+p_{12}(1-q^{-3})(1-q^{-2})(1+q)g_1f_1g_2^2f_2^2x_2= -p_{12}q^3(1-q^{-3})(1-q^{-2})(1+q)(1+q-q^{-1}-q^{-2}+q^{-1}-1)g_1f_1g_2^2f_2^2x_2=-p_{12}q^4(1+q)(1-q^{-2})(1-q^{-3})^2g_1f_1g_2^2f_2^2x_2
$,

\item [(A42)]
$[[[[x_1,x_2],x_2],x_2],[x_2^-,[x_2^-,[x_2^-,x_1^-]]]]=[[[x_1,x_2],x_2],[x_2,[x_2^-,[x_2^-,[x_2^-,x_1^-]]]]]+p_{22}^3p{12}[[[[x_1,x_2],x_2],[x_2^-,[x_2^-,[x_2^-,x_1^-]]]],x_2]=q_2(1-q^{-3})[[x_1,x_2],x_2][x_2^-,[x_2^-,x1^-]]$ $g_2f_2-p_{12}p_{22}^2p_{21}^3p_{22}^6p_{11}p_{12}^2q_2(1-q^{-3})[x_2^-,[x_2^-,x_1^-]]g_2f_2[[x_1,x_2],x_2]+p_{12}q^3p_{21}(1-q^{-3})(1-q^{-2})(1+q)(1-q^3)x_2^-x_2-p_{12}p_{22}p_{22}^{-3}p_{12}^{-1}p_{12}q^3p_{21}(1-q^{-3})(1-q^{-2})(1+q)(1-q^3)x_2x_2^-=q^2(1-q^{-3})(1-q^{-3})(1-q^{-2})(1+q)(1-g_1f_1g_2^2f_2^2)+p_{21}^2p_{22}^4p_{11}p_{12}^2[x_2^-,[x_2^-,x_1^-]][[x_1,x_2],x_2]g_2f_2-q^4(1-q^{-3})[x_2^-,[x_2^-,x_1^-]]$ $g_2f_2[[x_1,x_2],x_2]+(1-q^{-3})(1-q^{-2})(1+q)(1-q^3)x_2^-x_2-q^{-1}(1-q^{-3})(1-q^{-2})(1+q)(1-q^3)(1-g_2f_2+p_{22}x_2^-x_2)=q^2(1-q^{-3})^2(1-q^{-2})(1+q)(g_2f_2-g_1f_1g_2^3f_2^3)+q^4(1-q^{-3})[x_2^-,[x_2^-,x_1^-]]g_2f_2[[x_1,x_2],x_2]-q^4(1-q^{-3})[x_2^-[x_2^-,x_1^-]]g_2f_2[[x_1,x_2],x_2]-q^{-1}(1-q^{-3})(1-q^{-2})(1+q)(1-q^3)(1-g_2f_2)=q^2(1+q)(1-q^{-2})(1-q^{-3})^2(1-g_1f_1g_2^3f_2^3)$,

\item [(A43)]
$[[x_2,[x_2,[x_2,x_1]]],x_2^-]=[x_2[x_2,[x_2,x_1]],x_2^-]-p_{22}^2p_{21}[[x_2,[x_2,x_1]]x_2,x_2^-]=p_{22}^2p_{21}[x_2,x_2^-][x_2,[x_2,x_1]]+x_2[[x_2,[x_2,x_1]],x_2^-]-p_{22}^3p_{21}[[x_2,[x_2,x_1]],x_2^-]x_2-p_{22}^2p_{21}[x_2,[x_2,x_1]][x_2,x_2^-]=p_{21}q^2(1-g_2f_2)[x_2,[x_2,x_1]]+p_{21}(1-q+q^{-1}-q^2)x_2g_2f_2[x_2,x_1]-q^3p_{21}^2(1-q+q^{-1}-q^2)g_2f_2[x_2,x_1]x_2-p_{21}q^2[x_2,[x_2,x_1]](1-g_2f_2)=-p_{21}q^2g_2f_2[x_2,[x_2,x_1]]+p_{21}q^3g_2f_2[x_2,[x_2,x_1]]+p_{21}q^2(1-q+q^{-1}-q^2)g_2f_2x_2[x_2,x_1]-p_{21}^2q^3(1-q+q^{-1}-q^2)g_2f_2[x_2,x_1]x_2=p_{21}q(1-q^{-3})g_2f_2[x_2,[x_2,x_1]]$,

\item [(A44)]
$[[x_2,[x_2,[x_2,x_1]]],[[x_1^-,x_2^-],x_2^-]]=[[[x_2,[x_2,[x_2,x_1]]],[x_1^-,x_2^-]],x_2^-]+p_{12}^3p_{11}$ $p_{22}^3p_{21}[[x_1^-,x_2^-],[[x_2,[x_2,[x_2,x_1]]],x_2^-]]=q^2(1-q^{-3})^2(1-q^{-2})x_2^2x_2^--p_{22}^3p_{21}p_{21}^{-1}p_22^{-1}q^2(1-q^{-3})^2(1-q^{-2})x_2^-x_2^2+p_{12}^2q^3p_{21}q(1-q^3)[x_1^-,x_2^-]g_2f_2[x_2,[x_2,x_1]]-p_{12}^{-3}p_{11}^{-1}p_{22}^{-3}p_21^{-1}p_21^{-1}p_{22}^{-1}p_{12}^2q^3p_{21}q(1-q^3)g_2f_2[x_2,[x_2,x_1]][x_1,x_2^-]=q^2(1-q^{-3})^2(1-q^{-2})((1+q)x_2-q(1+q)g_2f_2x_2+q^2x_2^-x_2^2)-q^4(1-q^{-3})^2(1-q^{-2})x_2^-x_2^2+p_{12}q(1-q^3)[x_1^-,x_2^-]g_2f_2[x_2,[x_2,x_1]]-(1-q^3)g_2f_2((1-q^{-3}(1-q^{-2}(1+q)x_2+p_{12}^2p_{22}^2p_{11}p_{21}[x_1^-,x_2^-][x_2,[x_2,x_1]]=q^2(1-q^{-3})^2(1-q^{-2})(1+q)x_2-q^3(1-q^{-3})^2(1-q^{-2})(1+q)g_2f_2x_2+q^3(1-q^{-3})^2(1-q^{-2})(1+q)g_2f_2x_2=q^2(1+q)(1-q^{-2})(1-q^{-3})x_2$,

\item [(A45)]
$[[x_2,[x_2,[x_2,x_1]]],[[[x_1^-,x_2^-],x_2^-],x_2^-]]=[[[x_2,[x_2,[x_2,x_1]]],[[x_1^-,x_2^-],x_2^-]],x_2^-]+p_{12}^3p_{11}p_{22}^6p_{21}^2[[[x_1^-,x_2^-],x_2^-],[[x_2,[x_2,[x_2,x_1]]],x_2^-]]=q^2(1-q^{-3})^2(1-q^{-2})(1+q)x_2x_2^--p_{22}^3p_{21}p_{21}^{-1}p_{22}^{-2}q^2(1-q^{-3})^2(1-q^{-2})(1+q)x_2^-x_2+p_{12}q^3p_{21}q(1-q^3)[[x_1^-,x_2^-],x_2^-]g_2f_2[x_2,[x_2,x_1]]-p_{21}^{-1}p_{22}^{-2}p_{12}^{-3}p_{11}^{-1}p_{22}^{-6}p_{21}^{-2}q(1-q^3)g_2f_2((1-q^{-3})(1-q^{-2})(1+q)(1-g_1f_1g_2^2f_2^2)+p_{12}^2p_{11}p_{22}^4p_{21}^2[[x_1^-,x_2^-],x_2^-][x_2,[x_2,x_1]])=q^2(1-q^{-3})^2(1-q^{-2}(1+q)(1-g_2f_2)-q^{-1}(1-q^3)(1-q^{-3})(1-q^{-2})(1+q)(g_2f_2-g_1f_1g_2^3f_2^3)=q^2(1+q)(1-q^{-2})(1-q^{-3})^2(1-g_1f_1g_2^3f_2^3)$,

\item [(A46)]
$[[[x_1,x_2],[x_2,[x_2,x_1]]],x_2^-]=[[x_1,x_2],[[x_2,[x_2,x_1]],x_2^-]]+p_{21}p_{22}^2[[[x_1,x_2],x_2^-],$ $[x_2,[x_2,x_1]]]=p_{21}(1-q+q^{-1}-q^2)[x_1,x_2]g_2f_2[x_2,x_1]-p_{12}^2p_{11}p_{22}^3p_{21}^3(1-q+q^{-1}-q^2)g_2f_2[x_2,x_1][x_1,x_2]+p_{21}q^2(1-q^{-3})x_1[x_2,[x_2,x_1]]-p_{12}^2p_{11}p_{22}^2p_{21}p_{22}^{-2}p_{21}^{-1}p_{21}q^2(1-q^{-3})[x_2,[x_2,x_1]]x_1=p_{21}q^{-1}(1-q+q^{-1}-q^2)g_2f_2(q^3p_{12}(1-q^{-3})x_2x_1-p_{12}q^3[x_2,x_1])[x_2,x_1]-p_{21}(1-q+q^{-1}-q^2)g_2f_2[x_2,x_1]q^3p_{12}(1-q^{-3})x_2x_1-p_{12}q^3[x_2,x_1])+p_{21}q^2(1-q^{-3})(p_{12}q[x_2,x_1]^2+p_{12}^2q^3[x_2,[x_2,x_1]]x_1)-p_{12}q^2(1-q^{-3})[x_2,[x_2,x_1]]x_1=(1-q^{-3})[x_2,x_1]^2+p_{12}q^3(1+q)(1-q^{-2})(1-q^{-3})^2g_2f_2x_2[x_2,x_1]x_1-p_{12}q^3(1+q)(1-q^{-2})(1-q^{-3})g_2f_2[x_2,[x_2,x_1]]x_1+(1+q)^2(1-q^{-2})g_2f_2[x_2,x_1]^2$,

\item [(A47)]
$[[[x_1,x_2],[[x_1,x_2],x_2]],x_2^-]=[[x_1,x_2],[[[x_1,x_2],x_2],x_2^-]]+p_{21}p_{22}^2[[[x_1,x_2],x_2^-],$ $[[x_1,x_2],x_2]]=(1+q-q^{-1}-q^{-2})[x_1,x_2]^2-p_{11}p_{12}^2p_{21}p_{22}^2p_{21}p_{22}(1+q-q^{-1}-q^{-2})[x_1,x_2]^2+p_{21}q^2(1-q^{-3}x_1[[x_1,x_2],x_2]-p_{11}p_{12}^2p_{21}p_{22}^2p_{21}^{-1}p_{22}^{-2}p_{21}q^2(1-q^{-3})[[x_1,x_2],x_2]x_1=-p_{12}q^2(1-q^{-3})[[x_1,x_2],x_2]x_1+p_{21}q^2(1-q^{-3})(p_{12}q(1-q^{-2})[x_1,x_2]^2+p_{12}^2q^3[[x_1,x_2],x_2]x_1)=(1-q^{-2})(1-q^{-3})[x_1,x_2]^2$,

\item [(A48)]
$[[[x_1,x_2],[[x_1,x_2],x_2]],x_1^-]=[[x_1,x_2],[[[x_1,x_2],x_2],x_1^-]]+p_{11}p_{12}[[[x_1,x_2],x_1^-],$ $[[x_1,x_2],x_2]]=-p_{12}^2(1-q^{-3})(1-q^{-2})[x_1,x_2]g_1f_1x_2^2+p_{11}p_{12}^2p_{21}p_{22}^2p_{11}p_{12}p_{12}^2(1-q^{-3})(1-q^{-2})g_1f_1x_2^2[x_1,x_2]-p_{12}^2p_{11}(1-q^{-3})g_1f_1x_2[[x_1,x_2],x_2]+p_{11}p_{12}^2p_{21}p_{22}^2$ $p_{11}^{-1}p_{12}^{-2}p_{12}^{-2}p_{11}(1-q^{-3}[[x_1,x_2],x_2]g_1f_1x_2=-p_{12}^2q^3(1-q^{-3})(1-q^{-2})g_1f_1([[[x_1,x_2],x_2],x_2]+p_{12}q(1+q)x_2[[x_1,x_2],x_2]+p_{12}^2q^2x_2^2[x_1,x_2])+p_{12}^4q^5(1-q^{-3})(1-q^{-2})g_1f_1x_2^2[x_1,x_2]-p_{12}^2q^3(1-q^{-3})g_1f_1x_2[[x_1,x_2],x_2]+p_{12}q^2(1-q^{-3})g_1f_1([[[x_1,x_2],x_2],x_2]+p_{12}q^2x_2[[x_1,x_2],x_2])=-p_{12}^2q^3(1-q^{-3})g_1f_1[[[x_1,x_2],x_2],x_2]-p_{12}^3q^4(1+q)(1-q^{-2})(1-q^{-3})g_1f_1x_2[[x_1,x_2],x_2]$,

\item [(A49)]
$[[[x_1,x_2],[[x_1,x_2],x_2]],[x_2^-,x_1^-]]=[[x_1,x_2],[[[x_1,x_2],[x_2^-,x_1^-]]]+p_{21}p_{22}p_{11}p_{12}$ $[[[x_1,x_2],[x_2^-,x_1^-]],[[x_1,x_2],x_2]]=-p_{12}q(1-q^{-3})(1-q^{-2})(1+q)[x_1,x_2]g_1f_1g_2f_2x_2+p_{11}p_{12}^2p_{21}p_{22}^2p_{21}p_{22}p_{11}p_{12}p_{12}q(1-q^{-3})(1-q^{-2})(1+q)g_1f_1g_2f_2x_2[x_1,x_2]+q(1-q^{-3})(1-g_1f_1g_2f_2)[[x_1,x_2],x_2]-p_{11}p_{12}^2p_{21}p_{22}^2p_{21}^{-1}p_{22}^{-2}p_{11}^{-1}p_{12}^{-2}q(1-q^{-3}[[x_1,x_2],x_2](1-g_1f_1g_2f_2)=-p_{12}q^3(1-q^{-3})(1-q^{-2})(1+q)g_1f_1g_2f_2([[x_1,x_2],x_2]-p_{12}qx_2[x_1,x_2])+p_{12}^2q^4(1-q^{-3})(1-q^{-2})(1+q)g_1f_1g_2f_2x_2[x_1,x_2]+q(1-q^{-3})(1-g_1f_1g_2f_2)[[x_1,x_2],x_2]-q(1-q^{-3})[[x_1,x_2],x_2]+q(1-q^{-3})qg_1f_1g_2f_2[[x_1,x_2],x_2]=-p_{12}q^4(1-q^{-3})^2g_1f_1g_2f_2[[x_1,x_2],x_2]$,

\item [(A50)]
$[[[x_1,x_2],[[x_1,x_2],x_2]],[x_2^-,[x_2^-,x_1^-]]]=[[x_1,x_2],[[[x_1,x_2],x_2],[x_2^-,[x_2^-,x_1^-]]]]+p_{21}^2p_{22}^4p_{11}p_{12}^2[[[x_1,x_2],[x_2^-,[x_2^-,x_1^-]]],[[x_1,x_2],x_2]]=(1-q^{-3})(1-q^{-2})(1+q)[x_1,x_2](1-g_1f_1g_2^2f_2^2)-p_{11}p_{12}^2p_{21}p_{22}^2p_{21}^2p_{22}^2p_{11}p_{12}(1-q^{-3})(1-q^{-2})(1+q)(1-g_1f_1g_2^2f_2^2)[x_1,x_2]+qp_{21}(1-q^{-3})(1-q^2)(1+q)x_2^-[[x_1,x_2],x_2]-p_{11}p_{12}^2p_{21}p_{22}^2p_21^{-2}p_{22}^{-4}p_{11}^{-1}p_{12}^{-2}qp_{21}(1-q^{-3})(1-q^2)(1+q)((1+q-q^{-1}-q^{-2})[x_1,x_2]+p_{21}q^2x_2^-[[x_1,x_2],x_2])=(1-q^{-3})(1-q^{-2})(1+q)(1-q)[x_1,x_2]+p_{21}q(1-q^{-3})(1-q^2)(1+q)x_2^-[[x_1,x_2],x_2]-q^{-1}(1-q^{-3})(1-q^2)(1+q)(1+q-q^{-1}-q^{-2})[x_1,x_2]+p_{21}q(1-q^{-3})(1-q^2)(1+q)x_2^-[[x_1,x_2],x_2]=q^2(1+q)(1-q^{-2})(1-q^{-3})^2[x_1,x_2]$.

\end{description}

\section*{Acknowledgments}

I would like to thank professor Kharchenko for all his attention
during my research period in Mexico. I also thank him for
proposing this theme as part of my thesis and for all his suggestions to
this article.

\end{document}